\documentclass[12pt]{article}
\usepackage{amssymb,amsmath,amsthm,citesort}
\numberwithin{equation}{section}

\makeatletter
\renewcommand{\@makefnmark}{\hbox{\mathsurround=0pt $^\ast$}}
\renewcommand{\@biblabel}[1]{#1.\hfill}
\makeatother

\renewcommand*{\leq}{\leqslant} \renewcommand*{\geq}{\geqslant}

\newcommand*{\vare}{\varepsilon} \newcommand*{\vart}{\vartheta}
\newcommand*{\vark}{\varkappa} \newcommand*{\varp}{\varphi}

\newcommand*{\Gbig}{\Gamma_{\mathrm{big}}}

\newcommand*{\mC}{\mathbb C} \newcommand*{\mN}{\mathbb N}
\newcommand*{\mR}{\mathbb R} \newcommand*{\mT}{\mathbb T}
\newcommand*{\mZ}{\mathbb Z}

\newcommand*{\cB}{\mathcal B} \newcommand*{\cC}{\mathcal C}
\newcommand*{\cD}{\mathcal D} \newcommand*{\cG}{\mathcal G}
\newcommand*{\cJ}{\mathcal J} \newcommand*{\cK}{\mathcal K}
\newcommand*{\cL}{\mathcal L} \newcommand*{\cM}{\mathcal M}
\newcommand*{\cO}{\mathcal O} \newcommand*{\cQ}{\mathcal Q}
\newcommand*{\cR}{\mathcal R} \newcommand*{\cS}{\mathcal S}
\newcommand*{\cT}{\mathcal T} \newcommand*{\cU}{\mathcal U}

\newcommand*{\fB}{\mathfrak B} \newcommand*{\fC}{\mathfrak C}
\newcommand*{\fD}{\mathfrak D} \newcommand*{\fF}{\mathfrak F}
\newcommand*{\fG}{\mathfrak G} \newcommand*{\fM}{\mathfrak M}
\newcommand*{\fO}{\mathfrak O} \newcommand*{\fQ}{\mathfrak Q}
\newcommand*{\fS}{\mathfrak S} \newcommand*{\fU}{\mathfrak U}
\newcommand*{\fW}{\mathfrak W}

\newcommand*{\fa}{\mathfrak a} \newcommand*{\fb}{\mathfrak b}
\newcommand*{\fw}{\mathfrak w} \newcommand*{\fz}{\mathfrak z}

\newcommand*{\bx}{\overline{x}} \newcommand*{\by}{\overline{y}}
\newcommand*{\bz}{\overline{z}}

\newcommand*{\bsigma}{\overline{\sigma}}

\newcommand*{\const}{\mathrm{const}} \newcommand*{\id}{\mathrm{id}}
\newcommand*{\Ad}{\mathrm{Ad}}

\newcommand*{\gl}{\mathfrak{gl}} \newcommand*{\GL}{\mathrm{GL}}

\newcommand*{\codim}{\mathop{\mathrm{codim}}\nolimits}
\newcommand*{\diag}{\mathop{\mathrm{diag}}\nolimits}
\newcommand*{\Fix}{\mathop{\mathrm{Fix}}\nolimits}
\newcommand*{\rIm}{\mathop{\mathrm{Im}}\nolimits}

\newtheorem{lemma}{Lemma}
\newtheorem{thrm}{Theorem}

\theoremstyle{definition}
\newtheorem*{defn}{Definition}
\newtheorem*{remark}{Remark}

\allowdisplaybreaks

\begin{document}
\baselineskip=21pt
\parindent=0mm

\begin{center}
{\Large\bfseries Whitney Smooth Families of Invariant Tori \\
within the Reversible Context~2 of KAM Theory}

\bigskip
\baselineskip=16pt

{\large\bfseries Mikhail B.~Sevryuk\footnote{E-mails: \texttt{sevryuk@mccme.ru, 2421584@mail.ru}}}

\medskip

{\small\itshape V.~L.~Talroze Institute of Energy Problems of Chemical Physics of the Russia Academy of Sciences, \\
Leninskii prospect~38, Building~2, Moscow 119334, Russia}
\end{center}

\bigskip
\baselineskip=16pt

\textbf{Abstract}---%
We prove a general theorem on the persistence of Whitney $C^\infty$-smooth families of invariant tori in the reversible context~2 of KAM theory. This context refers to the situation where $\dim\Fix G < (\codim\cT)/2$ where $\Fix G$ is the fixed point manifold of the reversing involution $G$ and $\cT$ is the invariant torus in question. Our result is obtained as a corollary of the theorem by H.~W.~Broer, M.-C.~Ciocci, H.~Han{\ss}mann, and A.~Vanderbauwhede of 2009 concerning quasi-periodic stability of invariant tori with singular ``normal'' matrices in reversible systems.

\bigskip

MSC2010 numbers: \texttt{70K43, 70H33}

\bigskip


\bigskip

\begin{flushright}\itshape
To the blessed memory of Nikola\u{\i} Nikolaevich Nekhoroshev, \\
an outstanding mathematician and a wonderful personality
\end{flushright}

\parindent=5mm

\section{Introduction}
\label{intro}

\subsection{Reversible Systems and Their Invariant Tori}
\label{revsyst}

KAM theory was founded in the fifties and sixties of the 20th century by A.~N.~Kolmogorov, V.~I.~Arnold, and J.~Moser as the theory of quasi-periodic motions in non-integrable Hamiltonian systems (Hamiltonian flows and symplectic diffeomorphisms). However, it was soon realized that almost all the concepts and results of this theory can be carried over to other classes of dynamical systems, in particular, to reversible, volume preserving, and general (dissipative) systems. From this viewpoint, one sometimes speaks of the Hamiltonian, reversible, volume preserving, and dissipative \emph{contexts} of KAM theory. On the other hand, the properties of invariant tori filled up densely by quasi-periodic motions (conditionally periodic motions with incommensurable frequencies) depend strongly on the phase space structures the system in question is assumed to preserve. For recent general reviews of KAM theory the reader is referred to e.g.\ the tutorial \cite{dL01}, the monograph \cite[\S~6.3]{AKN06}, and the survey \cite{BS10}. As a semi-popular introduction to the theory, the book \cite{D14} is highly recommended. Some very recent achievements are expounded in the monograph \cite{HCLMF16}. There are many works on KAM theory where several different contexts are considered and the corresponding results are presented in parallel, see e.g.\ \cite{BHT90,S95JDCS,BHS96Gro,BHS96LNM,S06,S07Stek,S07DCDS,BHN07}.

Of all non-Hamiltonian versions of KAM theory, the context most similar to the Hamiltonian one is the reversible context. Recall that a dynamical system with a phase space $\cM$ equipped with a smooth \emph{involution} $G:\cM\to\cM$ (a mapping whose square is the identical transformation) is said to be \emph{reversible} with respect to $G$ (or $G$-reversible) if this system is invariant under the transformation $(w,t)\mapsto(Gw,-t)$ where $w\in\cM$ and $t$ is the time. For instance, the autonomous flow on $\cM$ afforded by a vector field $V$ is reversible with respect to $G$ if $\cD G(V\circ G)=-V$ where $\cD G$ denotes the differential of $G$. If this is the case, the field $V$ is also said to be $G$-reversible. General surveys of the theory of finite dimensional reversible systems with many examples and extensive bibliographies are given in the papers \cite{RQ92,LR98}.

The reversible KAM theory was founded in the mid sixties by J.~Moser, Yu.~N.~Bibikov, and V.~A.~Pliss \cite{M65,M66,M67,BP67}. By now, the literature on the reversible context of KAM theory is enormous, see e.g.\ the works \cite{S95JDCS,BHS96Gro,BHS96LNM,S06,S07Stek,S07DCDS,BHN07} cited above and the papers \cite{BH95,S95Cha,S98,L01,W01,BCHV09,WX09,H11RCD,H11IM}. In particular, the so-called parametric approach to KAM theory developed for the Hamiltonian, volume preserving, and dissipative contexts in \cite{BHT90} was carried over to reversible systems in \cite{BH95}. The paper \cite{S98} contains a brief review of the reversible KAM theory as it stood in 1997. Many references on KAM theory for reversible systems are presented in \cite{S11,S12MMJ,S12IM}. The very recent studies are exemplified by the papers \cite{WXZ15,WXZ16}. The reversible KAM theory for infinite dimensional systems is also known. Some partial results in this field were obtained as early as 1990 \cite{IL90}, the work \cite{CY07} is regarded as very important, and the recent contributions are exemplified by \cite{FP15}. However, in the present paper, we will always assume the phase space $\cM$ of a reversible system to be finite dimensional and connected.

If a set $\cK\subset\cM$ is invariant under a $G$-reversible system in $\cM$, so is the set $G(\cK)$. Nevertheless, in the case where $G(\cK)$ and $\cK$ do not coincide, the dynamical characteristics of these two invariant sets are opposite (for instance, one of them is an attractor whereas the other set is a repeller), and these sets exhibit \emph{no} special dynamical features compared with invariant sets of dissipative systems. Therefore, in the theory of reversible systems, one usually considers only sets invariant under \emph{both} the system itself and the reversing involution. In the sequel, while speaking of invariant tori of reversible systems, we will always comply with this rule.

Near a fixed point, any involution $G:\cM\to\cM$ is linear in suitable local coordinates (this is a very particular case of the Bochner theorem \cite{B72,MZ74}), and the set $\Fix G$ of fixed points of $G$ is a submanifold of $\cM$ of the same smoothness class as $G$ itself. This submanifold can well be empty or consist of several connected components of different dimensions even if $\cM$ is connected \cite{S11,B72,CF64} (in fact, the books \cite{B72,CF64} contain extensive information on the fixed point sets of involutions of various manifolds). However, for almost all the reversible systems encountered in practice, the fixed point manifold of the reversing involution $G$ is not empty and all the connected components of $\Fix G$ are of the same dimension, so that $\dim\Fix G$ is well defined \cite{RQ92,LR98}. In the present paper, we will only deal with autonomous reversible flows possessing this property.

The following lemma is very well known \cite{BHS96Gro,BHS96LNM,S12MMJ}.

\begin{lemma}\label{torus}
Let an $n$-torus $\cT\subset\cM$ be invariant under both a $G$-reversible flow on $\cM$ and the corresponding reversing involution $G$. If $\cT$ carries quasi-periodic motions then one can introduce a coordinate frame $\varp\in\mT^n=(\mR/2\pi\mZ)^n$ in $\cT$ such that the dynamics on $\cT$ takes the form $\dot{\varp}=\omega$ ($\omega\in\mR^n$ being the frequency vector of $\cT$) and the restriction of $G$ to $\cT$ takes the form $G|_{\cT}:\varp\mapsto-\varp$. Consequently, $\Fix G\cap\cT$ consists of $2^n$ isolated points $(\varp_1,\ldots,\varp_n)$ where each $\varp_j$ ($1\leq j\leq n$) is equal either to $0$ or to $\pi$, and $\dim\Fix G \leq \dim\cM-n=\codim\cT$.
\end{lemma}

Throughout the paper, $\cO_N(a)$ will denote an unspecified neighborhood of a point $a\in\mR^N$. For any variable $w\in\cO_N(0)$, we will write $O(w)$ instead of $O\bigl(|w|\bigr)$ and $O_l(w)$ instead of $O\bigl(|w|^l\bigr)$ for $l\geq 2$. Similarly, we will write $O(w,w',w'',\ldots)$ instead of $O\bigl(|w|+|w'|+|w''|+\cdots\bigr)$ and $O_l(w,w',w'',\ldots)$ instead of $O\bigl(|w|^l+|w'|^l+|w''|^l+\cdots\bigr)$ for $l\geq 2$. For integer vectors $k\in\mZ^N$, we will write $|k|=|k_1|+|k_2|+\cdots+|k_N|$. The angle brackets $\langle{\cdot},{\cdot}\rangle$ will denote the standard inner product of real vectors.

Instead of $\{0\}$ with $0\in\mR^N$, we will sometimes write $\{0\in\mR^N\}$. The space of $N_1\times N_2$ real matrices will be denoted by $\mR^{N_1\times N_2}$ (in particular, $\mR^{N\times N}=\gl(N,\mR)$). The block-diagonal matrix with blocks $M_1,M_2,\ldots$ will be denoted by $M_1\oplus M_2\oplus\cdots$. As usual, the symbol $\oplus$ will be also employed for the direct sum of subspaces that pairwise have only the zero vector in common.

Recall the concept of reducible invariant tori which is of principal importance in KAM theory \cite{AKN06,BHS96Gro,BHS96LNM,S07Stek,M66,M67,S95Cha,S98,S12MMJ,S12IM}.

\begin{defn}
Let an invariant $n$-torus $\cT$ of some flow on an $(n+N)$-di\-men\-sional manifold $\cM$ carry conditionally periodic motions with frequency vector $\omega\in\mR^n$. This torus is said to be \emph{reducible} (or \emph{Floquet}) if the variational equation along $\cT$ can be reduced to a form with constant coefficients, i.e., if in a neighborhood of $\cT$, there exist coordinates $\bigl(x\in\mT^n, \, X\in\cO_N(0)\bigr)$ in which the torus $\cT$ itself is given by the equation $\{X=0\}$ and the dynamical system takes the \emph{Floquet form} $\dot{x}=\omega+O(X)$, $\dot{X}=\Omega X+O_2(X)$ with an $x$-independent matrix $\Omega\in\gl(N,\mR)$. This matrix (not determined uniquely) is called the \emph{Floquet matrix} of the torus $\cT$, and its eigenvalues are called the \emph{Floquet exponents} of $\cT$. The coordinates $(x,X)$ are called the \emph{Floquet coordinates} for $\cT$.
\end{defn}

Another key concept in KAM theory is that of Diophantine approximations. In this paper, we will use the following definition.

\begin{defn}
Let the non-zero eigenvalues of a matrix $Q\in\gl(N,\mR)$ come in pairs $(\fa,-\fa)$ (this is definitely so if $Q$ anti-commutes with a non-singular matrix). Suppose that the matrix $Q$ possesses $\ell\geq 0$ pairs of non-zero purely imaginary eigenvalues
\begin{equation}
\pm i\beta_1,\ldots,\pm i\beta_\ell
\label{imaginary}
\end{equation}
and $\vark\geq 0$ quadruplets of eigenvalues outside $\mR\cup i\mR$:
\begin{equation}
\pm\alpha_1\pm i\beta_{\ell+1},\ldots,\pm\alpha_\vark\pm i\beta_{\ell+\vark}
\label{quadruplet}
\end{equation}
(one may assume that all the numbers $\alpha_j$ and $\beta_j$ are positive). Let also $\omega\in\mR^n$ be a vector and $\tau$, $\gamma$ be positive numbers. The pair $(\omega,Q)$ is said to be $(\tau,\gamma)$-\emph{Dio\-phan\-tine} if the inequality
\begin{equation}
\bigl| \langle k,\omega\rangle + \langle K,\beta\rangle \bigr| \geq \gamma|k|^{-\tau}
\label{Dioph}
\end{equation}
holds for all $k\in\mZ^n\setminus\{0\}$ and $K\in\mZ^{\ell+\vark}$ with $|K|\leq 2$.
\end{defn}

Let $n\geq 2$. It is well known that if $\omega_1,\ldots,\omega_n, \, \beta_1,\ldots,\beta_{\ell+\vark}$ are independent quantities and the point $(\omega,\beta)$ ranges in a bounded open domain in $\mR^{n+\ell+\vark}$, then for any fixed $\tau>n-1$ the infinite system of inequalities~\eqref{Dioph} determines a Cantor-like subset $\fD$ (a nowhere dense set of positive Lebesgue measure) of this domain (for $\gamma$ sufficiently small). To be more precise, the Lebesgue measure of the complement of $\fD$ is $O(\gamma)$ as $\gamma\to 0$ \cite{BHS96Gro,BHS96LNM,S06,S07Stek,S95Cha}. Note that if a pair $(\omega,Q)$ is $(\tau,\gamma)$-Dio\-phan\-tine then the vector $\omega$ is $(\tau,\gamma)$-Dio\-phan\-tine in the usual sense.

\subsection{Reversible Contexts~1 and~2}
\label{revcont}

An overwhelming majority of the works on the reversible KAM theory is devoted to the following problem. One considers systems of ordinary differential equations of the form
\begin{equation}
\dot{x}=H(Y,\lambda)+O(z), \quad \dot{Y}=O(z), \quad \dot{z}=Q(x,Y,\lambda)z+O_2(z)
\label{eq1}
\end{equation}
reversible with respect to the involution
\begin{equation}
G:(x,Y,z)\mapsto(-x,Y,Rz),
\label{rev1}
\end{equation}
where $x\in\mT^n$, $Y\in\mR^m$, $z\in\cO_{2p}(0)$ are the phase space variables, $\lambda\in\mR^s$ is an external parameter ($n$, $m$, $p$, $s$ being non-negative integers), $R\in\GL(2p,\mR)$ is an involutive matrix with eigenvalues $1$ and $-1$ of multiplicity $p$ each, $Q$ is a $2p\times 2p$ matrix-valued function, and $RQ(-x,Y,\lambda)R\equiv-Q(x,Y,\lambda)$. For each value of $\lambda$, system~\eqref{eq1} and involution~\eqref{rev1} admit the family $\{Y=\const, \, z=0\}$ of invariant $n$-tori carrying conditionally periodic motions with frequency vectors $H(Y,\lambda)$, and one is looking for invariant $n$-tori close to $\{Y=\const, \, z=0\}$ in small $G$-reversible perturbations of family~\eqref{eq1}.

In the case where $p\geq 1$ and the matrices $Q=Q(Y,\lambda)\in\gl(2p,\mR)$ are $x$-independent and possess $\ell$ pairs of non-zero purely imaginary eigenvalues ($1\leq\ell\leq p$) for some $Y$ and $\lambda$, one also examines invariant $(n+d)$-tori for $1\leq d\leq\ell$ near the surface $\{z=0\}$ in family~\eqref{eq1} itself as well as in its small $G$-reversible perturbations (the so-called excitation of elliptic normal modes, see \cite{BHS96LNM,S95Cha,S98} and references therein). It is also possible to study bifurcations of invariant $n$-tori into invariant $(n+d)$-tori ($1\leq d\leq p$) in $G$-reversible systems close to Eq.~\eqref{eq1} (see \cite{H11RCD,H11IM} and references therein).

Now note that the phase space codimension of an invariant $(n+d)$-torus $\cT$ here ($0\leq d\leq p$) is equal to $m+2p-d$ while $\dim\Fix G=m+p$. Thus, for involution~\eqref{rev1} and $(n+d)$-tori in question one has
\begin{equation}
\frac{\codim\cT}{2} \leq \dim\Fix G \leq \codim\cT
\label{context1}
\end{equation}
(this is tantamount to that $p+(m-d)/2\leq m+p\leq m+2p-d$). However, the two inequalities in Eq.~\eqref{context1} are of \emph{entirely different nature}. The right inequality follows from Lemma~\ref{torus} and holds for any invariant torus $\cT$ (carrying quasi-periodic motions) of any flow reversible with respect to any involution $G$. On the other hand, the left inequality in Eq.~\eqref{context1} stems just from the special form of systems~\eqref{eq1} and involution~\eqref{rev1}.

Indeed, nothing prevents one from considering e.g.\ systems of the form
\begin{equation}
\dot{x}=\omega(\lambda)+O(y,z), \quad \dot{y}=\sigma(\lambda)+O(y,z), \quad \dot{z}=Q(x,\lambda)z+O_2(y,z)
\label{eq2}
\end{equation}
reversible with respect to the involution
\begin{equation}
G:(x,y,z)\mapsto(-x,-y,Rz),
\label{rev2}
\end{equation}
where $x\in\mT^n$, $y\in\cO_m(0)$, $z\in\cO_{2p}(0)$ are the phase space variables, $\lambda\in\mR^s$ is an external parameter ($n$, $p$, $s$ being non-negative integers while $m$ being a positive integer), $R\in\GL(2p,\mR)$ is an involutive matrix with eigenvalues $1$ and $-1$ of multiplicity $p$ each, $Q$ is a $2p\times 2p$ matrix-valued function, and $RQ(-x,\lambda)R\equiv-Q(x,\lambda)$. The reversibility with respect to involution~\eqref{rev2} does not preclude a drift along the variable $y$, but such a drift can be overcome by a multidimensional external parameter. To be more precise, for each value of $\lambda$ such that $\sigma(\lambda)=0$ (if $s\geq m$ then $\sigma^{-1}(0)$ is generically an $(s-m)$-di\-men\-sional surface in $\mR^s$), system~\eqref{eq2} and involution~\eqref{rev2} admit the invariant $n$-torus $\{y=0, \, z=0\}$ carrying conditionally periodic motions with frequency vector $\omega(\lambda)$. One may then try to construct invariant $n$-tori $\cT$ close to $\{y=0, \, z=0\}$ in small $G$-reversible perturbations of family~\eqref{eq2}. The phase space codimension of $\cT$ is equal to $m+2p$ while $\dim\Fix G=p$. So, for involution~\eqref{rev2} and $n$-tori in question one gets
\begin{equation}
\dim\Fix G < \frac{\codim\cT}{2}
\label{context2}
\end{equation}
($p<p+m/2$ because $m\geq 1$). We use the upper case letter $Y$ in Eqs.~\eqref{eq1}--\eqref{rev1} and the lower case letter $y$ in Eqs.~\eqref{eq2}--\eqref{rev2} to emphasize the difference between the involutions~\eqref{rev1} and~\eqref{rev2}.

\begin{defn}
Let $\cT$ be a torus invariant under a system reversible with respect to an involution $G$. The situation where the inequalities~\eqref{context1} hold is called \emph{the reversible context~1} while the opposite situation where the inequality~\eqref{context2} holds is called \emph{the reversible context~2}.
\end{defn}

This definition was introduced in the works \cite{BHS96Gro,BHS96LNM}. The differences between the reversible contexts~1 and~2 were discussed in detail in the paper \cite{S11}. By now, KAM theory for the reversible context~1 is nearly as developed as the Hamiltonian KAM theory. The task of initiating KAM theory for the reversible context~2 was stated as Problem~9 in the note \cite{S08}. Up to now, the only works where the reversible KAM theory in context~2 is dealt with have been the papers \cite{S11,S12MMJ,S12IM}. In \cite{S11}, we examined the so-called extreme reversible context~2 where $\dim\Fix G=0$ (i.e., $\Fix G$ is a finite collection of isolated points). In \cite{S12MMJ}, the general reversible context~2 for an arbitrary dimension of $\Fix G$ satisfying Eq.~\eqref{context2} was studied. Finally, the paper \cite{S12IM} considered non-autonomous systems (depending on time quasi-periodically) within the reversible context~2.

The principal technical tool in the papers \cite{S11,S12MMJ} was Moser's modifying terms theory \cite{M67}. In particular, in \cite{S12MMJ}, we obtained $(s-n-m-\ell-\vark)$-parameter \emph{analytic} families of reducible invariant $n$-tori $\cT$ ($n\geq 1$), where $s\geq n+m+\ell+\vark$ is the number of external parameters in the analytic $G$-reversible system in question, $m=\codim\cT-2\dim\Fix G>0$, and the Floquet matrix $\Omega$ of each torus $\cT$ possesses $\ell\geq 0$ pairs of non-zero purely imaginary eigenvalues~\eqref{imaginary} and $\vark\geq 0$ quadruplets~\eqref{quadruplet} of eigenvalues outside $\mR\cup i\mR$. Moreover, the frequencies $\omega_1,\ldots,\omega_n$ of the tori $\cT$ and the positive imaginary parts $\beta_1,\ldots,\beta_{\ell+\vark}$ of the Floquet exponents of these tori are \emph{the same} for all the tori $\cT$ in the given family, and the pairs $(\omega,\Omega)$ are $(\tau,\gamma)$-Dio\-phan\-tine for some constants $\tau>n-1$ and $\gamma>0$.

However, this result cannot be regarded as quite satisfactory. By analogy with the four ``well-developed'' contexts of KAM theory (the Hamiltonian, volume preserving, and dissipative contexts as well as the reversible context~1) \cite{AKN06,BS10,BHT90,BHS96Gro,BHS96LNM,S06,S07Stek,BHN07,BH95,S95Cha,BCHV09}, one would expect that the $(s-n-m-\ell-\vark)$-parameter analytic families of reducible invariant $n$-tori constructed in \cite{S12MMJ} are organized into $(s-m)$-parameter \emph{Whitney smooth} families. The tori in such families have different frequencies $\omega$ and positive imaginary parts $\beta$ of the Floquet exponents and depend on $\omega$ and $\beta$ in a Whitney smooth way. In other words, the Floquet coordinates for the tori are defined a~priori for Diophantine pairs $(\omega,\beta)$ only, but these coordinates can be continued to smooth (say, $C^\infty$) functions $\fF$ defined in an \emph{open} domain of the appropriate Euclidean space. The subset $\fS$ of this domain corresponding to the pairs $(\omega,\beta)$ satisfying Eq.~\eqref{Dioph} (and, consequently, to the invariant tori) is a Cantor-like set. The values of $\fF$ in the complement of $\fS$ have no dynamical meaning. Unfortunately, the techniques of \cite{M67} do not enable one to obtain Whitney smooth families of invariant tori.

\begin{remark}
While speaking of the Hamiltonian context of KAM theory, one usually has in view only isotropic invariant tori. However, non-isotropic invariant tori of Hamiltonian systems are also known to be organized into Whitney smooth families, see e.g.\ \cite{KP01}.
\end{remark}

\subsection{Subject and Structure of the Paper}
\label{subject}

In the present paper, we construct Whitney $C^\infty$-smooth families of reducible invariant tori in the reversible context~2 of KAM theory for analytic systems that depend analytically on external parameters. Our main tool is the BCHV (Broer--Ciocci--Han{\ss}mann--Vanderbauwhede) theorem \cite{BCHV09} describing such families of invariant tori in a certain particular case of the reversible context~1 (namely, in the case of systems~\eqref{eq1} with \emph{singular} $x$-independent matrices $Q$). Note that the landmark article \cite{BCHV09} was published earlier than our first paper \cite{S11} on the reversible context~2. Whitney smooth families of invariant tori in the reversible context~2 could have been obtained from the very beginning\ldots

According to some general observations in KAM theory \cite{P00,W03,W10}, invariant tori in analytic (and even Gevrey regular) dynamical systems are expected to be organized into families that are Gevrey regular in the sense of Whitney, not just $C^\infty$. Gevrey regularity of the families of invariant tori for some analytic flows in the reversible context~1 was established in \cite{WX09}. There is no doubt that the Whitney $C^\infty$-smooth families of invariant tori constructed in the present paper are also Gevrey regular, but we do not prove this.

In fact, the present paper develops the parametric KAM theory in the reversible context~2 for multiple non-zero Floquet exponents. For the Hamiltonian context, volume preserving context, dissipative context, and reversible context~1, the parametric framework of KAM theory was worked out in \cite{BHT90,BH95} for simple non-zero Floquet exponents and in \cite{BHN07,BCHV09} for multiple non-zero Floquet exponents. The author also expects that Wagener's general parametrized version of Moser's modifying terms theory \cite{W10} can be carried over to the reversible context~2 (according to H.~Han{\ss}mann \cite{H11RCD,H11IM}, this is so for the reversible context~1).

\begin{remark}
As was pointed out in \cite{S11}, the reversible context~2 appears naturally in the ``conventional'' setup of the involution $G$~\eqref{rev1} and $G$-reversible systems~\eqref{eq1} if one examines the break-up of \emph{resonant} unperturbed invariant $n$-tori $\{Y=\const, \, z=0\}$. Indeed, suppose that a resonant unperturbed invariant $n$-torus gives rise to a finite collection of perturbed invariant $r$-tori carrying quasi-periodic motions where $1\leq r\leq n-1$ (for Hamiltonian systems, this phenomenon has been explored very well by now, see short reviews in \cite{AKN06,BS10,BHS96LNM} and references therein). The phase space codimension of these $r$-tori is equal to $n+m+2p-r$, and these tori pertain to the reversible context~2 whenever
\[
\dim\Fix G=m+p < \frac{n+m+2p-r}{2}=\frac{n+m-r}{2}+p,
\]
i.e., $r<n-m$. However, to the best of the author's knowledge, the break-up of resonant unperturbed invariant tori in reversible systems (for $r\geq 2$) has been studied by now in the case of $p=0$, $m=n$ only \cite{L01,W01} where the inequality $r<n-m$ cannot be satisfied.
\end{remark}

The paper is organized as follows. Section~\ref{mainsec} contains the precise formulation of our main result. In Section~\ref{sourcesec}, we recall one of the particular cases of the BCHV theorem \cite{BCHV09} we need to prove our result. This proof is presented in Section~\ref{proof}. Section~\ref{Ruesssec} discusses the simplest application of our main result, namely, the analogue of the R\"ussmann nondegeneracy theorem for the reversible context~2. The remarks of Section~\ref{conclu} conclude the paper.

\section{The Main Result}
\label{mainsec}

One of the key concepts in the BCHV theory \cite{BCHV09} and in the present study is versal unfoldings with respect to actions of Lie groups. Let a Lie group $\fG$ act smoothly on a manifold $\fM$. For simplicity, we will assume $\fM$ and $\fG$ to be finite dimensional. Consider a smooth mapping $\fU:\cO_s(\lambda_\ast)\to\fM$ and the point $M_\ast=\fU(\lambda_\ast)$. The mapping $\fU$ can be regarded as an $s$-parameter \emph{unfolding} (or \emph{deformation}) of $M_\ast$.

\begin{defn}
The mapping $\fU$ is called a \emph{versal unfolding} (or \emph{versal deformation}) of the point $M_\ast\in\fM$ with respect to the given action of $\fG$ on $\fM$ if this mapping is transversal at $\lambda_\ast$ to the orbit $\fB(M_\ast)$ of $M_\ast$ under the action of $\fG$, i.e., if the tangent space $T_{M_\ast}\fM$ is spanned by $T_{M_\ast}\fB(M_\ast)$ and $\cD\fU(T_{\lambda_\ast}\mR^s)$ where $\cD\fU$ denotes the differential of $\fU$. The mapping $\fU$ is called a \emph{miniversal unfolding} (or \emph{miniversal deformation}) of $M_\ast$ if it is a versal unfolding of $M_\ast$ and the number $s$ of the parameters of this unfolding is minimum possible, i.e., $s$ is equal to the codimension of $\fB(M_\ast)$ in $\fM$ (so that $T_{M_\ast}\fM = T_{M_\ast}\fB(M_\ast) \oplus \cD\fU(T_{\lambda_\ast}\mR^s)$).
\end{defn}

The importance of versal unfoldings in singularity theory and related branches of mathematics is explained in detail in Arnold's works \cite{A71,A72,A88}. In particular, in those works, Arnold described versal unfoldings in the case where $\fM=\gl(N,\mC)$ is the space of all $N\times N$ complex matrices with the adjoint action (action by conjugations) of the group $\fG=\GL(N,\mC)$ of non-singular (invertible) $N\times N$ complex matrices: $\Ad_A(M)=AMA^{-1}$ for $M\in\gl(N,\mC)$, $A\in\GL(N,\mC)$. The unfolding parameters in Arnold's construction are also complex, i.e., instead of a real neighborhood $\cO_s(\lambda_\ast)$ of $\lambda_\ast\in\mR^s$, one considers a complex neighborhood of $\lambda_\ast\in\mC^s$ with a holomorphic mapping $\fU$.

In the sequel, we will need some information on versal unfoldings of real matrices anti-commuting with a fixed involutive matrix (of \emph{infinitesimally reversible} matrices in the terminology of \cite{S92}). Let $R\in\GL(N,\mR)$ be a fixed involutive matrix, $\fM=\gl_{-R}$ be the space of $N\times N$ real matrices anti-commuting with $R$, and $\fG=\GL_{+R}$ be the group of non-singular $N\times N$ real matrices commuting with $R$. The normal forms and versal unfoldings of matrices in $\gl_{-R}$ with respect to the adjoint action of $\GL_{+R}$ are found in \cite{S92,S93,I96}, some particular cases are also discussed in \cite{BCHV09}. An extensive bibliography on versal unfoldings for various matrix spaces is compiled in the book \cite[pp.~226--231]{A04} (to be more precise, in the comment to Problem \mbox{1970-1} in \cite{A04}).

Now we are in the position to present the main result of this paper. Consider an $S$-parameter analytic family of analytic systems
\begin{equation}
\begin{aligned}
\dot{x} &= \omega(\lambda)+\widetilde{\xi}(y,z,\lambda)+\widetilde{f}(x,y,z,\lambda), \\
\dot{y} &= \sigma(\lambda)+\widetilde{\eta}(y,z,\lambda)+\widetilde{g}(x,y,z,\lambda), \\
\dot{z} &= \widetilde{Q}(\lambda)z+\widetilde{\zeta}(y,z,\lambda)+\widetilde{h}(x,y,z,\lambda)
\end{aligned}
\label{eqstart}
\end{equation}
reversible with respect to the involution $G:(x,y,z)\mapsto(-x,-y,Rz)$, cf.\ Eqs.~\eqref{eq2}--\eqref{rev2}. Here $x\in\mT^n$, $y\in\cO_m(0)$, $z\in\cO_{2p}(0)$ are the phase space variables, $\lambda\in\cO_S(\lambda_\ast)$ is an external parameter ($n$, $p$ being non-negative integers while $m$, $S$ being positive integers), $R\in\GL(2p,\mR)$ is an involutive matrix with eigenvalues $1$ and $-1$ of multiplicity $p$ each, $\widetilde{Q}$ is a $2p\times 2p$ matrix-valued function, $R\widetilde{Q}(\lambda)\equiv-\widetilde{Q}(\lambda)R$, $\widetilde{\xi}=O(y,z)$, $\widetilde{\eta}=O_2(y,z)$, $\widetilde{\zeta}=O_2(y,z)$, and $\widetilde{f}$, $\widetilde{g}$, $\widetilde{h}$ are small perturbation terms.

Our main assumption is that the mapping
\[
\lambda \mapsto \left( \omega(\lambda),\sigma(\lambda),\widetilde{Q}(\lambda) \right) \in \mR^n\times\mR^m\times\gl_{-R}
\]
is a \emph{versal unfolding} of the triple $\left( \omega(\lambda_\ast),\sigma(\lambda_\ast),\widetilde{Q}(\lambda_\ast) \right)$ with respect to the natural action of the group $\{\id\}\times\{\id\}\times\GL_{+R}$, where $\{\id\}$ is the trivial group consisting of the identity element only. Moreover, we suppose that $\sigma(\lambda_\ast)=0$ and $\det\widetilde{Q}(\lambda_\ast)\neq 0$. Introduce the notations $\omega_\ast=\omega(\lambda_\ast)$ and $s=S-n-m\geq 0$. One can regard $\omega(\lambda)$ and $\sigma(\lambda)$ as independent parameters. To be more precise, as a new coordinate frame in the parameter space, one can choose $(\omega,\sigma,\mu)$ where
\[
\omega\in\cO_n(\omega_\ast), \quad \sigma\in\cO_m(0), \quad \mu\in\cO_s(0),
\]
and the point $(\omega=\omega_\ast, \, \sigma=0, \, \mu=0)$ corresponds to $\lambda=\lambda_\ast$. The systems~\eqref{eqstart} take the form
\begin{align*}
\dot{x} &= \omega+\xi(y,z,\omega,\sigma,\mu)+f(x,y,z,\omega,\sigma,\mu), \\
\dot{y} &= \sigma+\eta(y,z,\omega,\sigma,\mu)+g(x,y,z,\omega,\sigma,\mu), \\
\dot{z} &= \widehat{Q}(\omega,\sigma,\mu)z+\widehat{\zeta}(y,z,\omega,\sigma,\mu)+h(x,y,z,\omega,\sigma,\mu),
\end{align*}
where $\xi=O(y,z)$, $\eta=O_2(y,z)$, $\widehat{\zeta}=O_2(y,z)$. The mapping $\mu\mapsto\widehat{Q}(\omega_\ast,0,\mu)\in\gl_{-R}$ is a versal unfolding of the non-singular matrix $\widehat{Q}(\omega_\ast,0,0)=\widetilde{Q}(\lambda_\ast)$ with respect to the adjoint action of the group $\GL_{+R}$.

Now set
\begin{align*}
Q(\omega,\mu) &= \widehat{Q}(\omega,0,\mu), \\
\zeta(y,z,\omega,\sigma,\mu) &= \widehat{\zeta}(y,z,\omega,\sigma,\mu) + \left[ \widehat{Q}(\omega,\sigma,\mu)-\widehat{Q}(\omega,0,\mu) \right]z.
\end{align*}
We arrive at an $(n+m+s)$-parameter family of $G$-reversible systems
\begin{equation}
\begin{aligned}
\dot{x} &= \omega+\xi(y,z,\omega,\sigma,\mu)+f(x,y,z,\omega,\sigma,\mu), \\
\dot{y} &= \sigma+\eta(y,z,\omega,\sigma,\mu)+g(x,y,z,\omega,\sigma,\mu), \\
\dot{z} &= Q(\omega,\mu)z+\zeta(y,z,\omega,\sigma,\mu)+h(x,y,z,\omega,\sigma,\mu).
\end{aligned}
\label{eqfin}
\end{equation}
In this family, $\xi=O(y,z)$, $\eta=O_2(y,z)$, $\zeta=O_2(y,z,\sigma)$, and $RQ(\omega,\mu)\equiv-Q(\omega,\mu)R$. All the functions $Q$, $\xi$, $\eta$, $\zeta$, $f$, $g$, $h$ in Eq.~\eqref{eqfin} are analytic in all their arguments. The mapping $\mu\mapsto Q(\omega_\ast,\mu)\in\gl_{-R}$ is a \emph{versal unfolding} of the non-singular matrix $Q(\omega_\ast,0)=\widetilde{Q}(\lambda_\ast)$ with respect to the adjoint action of the group $\GL_{+R}$. In particular, this implies that $s \geq \codim\fB\bigl( Q(\omega_\ast,0) \bigr)$ where $\fB\bigl( Q(\omega_\ast,0) \bigr)$ is the orbit of $Q(\omega_\ast,0)$ under the adjoint action of $\GL_{+R}$. The minimum possible value of $\codim\fB\bigl( Q(\omega_\ast,0) \bigr)$ is equal to $p$ and attained in the case where to each eigenvalue of the matrix $Q(\omega_\ast,0)$, there corresponds a single Jordan block \cite{S92,S93,I96} (in \cite{A71}, such matrices were called Sylvester matrices).

Our main result will be formulated for the family~\eqref{eqfin} where the functions $Q$, $\xi$, $\eta$, $\zeta$ are fixed and the functions $f$, $g$, $h$ are small perturbation terms.

\begin{thrm}\label{main}
There exist a closed set $\Gamma\subset\mR^{n+s}$ that is diffeomorphic to an $(n+s)$-di\-men\-sional ball and contains the point $(\omega_\ast,0)$ in its interior, a number $\rho>0$, and a complex neighborhood $\cC \subset (\mC/2\pi\mZ)^n\times\mC^{n+2m+2p+s}$ of the set
\begin{equation}
\mT^n\times\{0\in\mR^m\}\times\{0\in\mR^{2p}\} \times \{\omega_\ast\}\times\{0\in\mR^m\}\times\{0\in\mR^s\}
\label{chord}
\end{equation}
with the following property. For any $\tau>n-1$, $\gamma>0$, $L\in\mN$, and $\vare>0$ there is $\delta>0$ such that the following holds. Suppose that the perturbation terms $f$, $g$, $h$ can be holomorphically continued to the neighborhood $\cC$ and $|f|<\delta$, $|g|<\delta$, $|h|<\delta$ in $\cC$. Then there exist mappings
\begin{equation}
\begin{aligned}
a &= a(\bx,\omega,\mu), &\quad a &: \mT^n\times\Gamma\to\mR^n, \\
b^0 &= b^0(\bx,\omega,\mu), &\quad b^0 &: \mT^n\times\Gamma\to\mR^m, \\
b^1 &= b^1(\bx,\omega,\mu), &\quad b^1 &: \mT^n\times\Gamma\to\gl(m,\mR), \\
b^2 &= b^2(\bx,\omega,\mu), &\quad b^2 &: \mT^n\times\Gamma\to\mR^{m\times 2p}, \\
c^0 &= c^0(\bx,\omega,\mu), &\quad c^0 &: \mT^n\times\Gamma\to\mR^{2p}, \\
c^1 &= c^1(\bx,\omega,\mu), &\quad c^1 &: \mT^n\times\Gamma\to\mR^{2p\times m}, \\
c^2 &= c^2(\bx,\omega,\mu), &\quad c^2 &: \mT^n\times\Gamma\to\gl(2p,\mR), \\
u &= u(\omega,\mu), &\quad u &: \Gamma\to\mR^n, \\
v &= v(\omega,\mu), &\quad v &: \Gamma\to\mR^m, \\
w &= w(\omega,\mu), &\quad w &: \Gamma\to\mR^s
\end{aligned}
\label{mapp}
\end{equation}
possessing the following properties. First, the mappings~\eqref{mapp} are analytic in $\bx$ and $C^\infty$-smooth in $\omega$ and $\mu$. Second, all the partial derivatives of each component of these mappings of any order from $0$ to $L$ are everywhere smaller than $\vare$ in absolute value. Third, for any $(\omega,\mu)\in\Gamma$, each component of $a$, $b^0$, $b^1$, $b^2$, $c^0$, $c^1$, and $c^2$ as a function of $\bx$ can be holomorphically continued to the strip
\begin{equation}
\cU_n(\rho) = \bigl\{ \bx\in(\mC/2\pi\mZ)^n \bigm| |\rIm\bx_j|<\rho, \; 1\leq j\leq n \bigr\}
\label{strip}
\end{equation}
and is smaller than $\vare$ in absolute value in this strip. Fourth, for any point $(\omega_0,\mu_0)\in\Gamma$, the nearly identical change of variables
\begin{equation}
\begin{aligned}
x &= \bx+a(\bx,\omega_0,\mu_0), \\
y &= \by+b^0(\bx,\omega_0,\mu_0)+b^1(\bx,\omega_0,\mu_0)\by+b^2(\bx,\omega_0,\mu_0)\bz, \\
z &= \bz+c^0(\bx,\omega_0,\mu_0)+c^1(\bx,\omega_0,\mu_0)\by+c^2(\bx,\omega_0,\mu_0)\bz
\end{aligned}
\label{change}
\end{equation}
with $\bx\in\mT^n$, $\by\in\cO_m(0)$, $\bz\in\cO_{2p}(0)$ commutes with $G$. Fifth, for any point $(\omega_0,\mu_0)\in\Gamma$ such that the pair $\bigl( \omega_0,Q(\omega_0,\mu_0) \bigr)$ is $(\tau,\gamma)$-Dio\-phan\-tine, the system~\eqref{eqfin} at the parameter values
\begin{equation}
\omega=\omega_0+u(\omega_0,\mu_0), \quad \sigma=v(\omega_0,\mu_0), \quad \mu=\mu_0+w(\omega_0,\mu_0)
\label{shift}
\end{equation}
after the coordinate change~\eqref{change} takes the form
\begin{equation}
\dot{\bx}=\omega_0+O(\by,\bz), \quad \dot{\by}=O_2(\by,\bz), \quad \dot{\bz}=Q(\omega_0,\mu_0)\bz+O_2(\by,\bz).
\label{goal}
\end{equation}
\end{thrm}

Consider a point $(\omega_0,\mu_0)\in\Gamma$ such that the pair $\bigl( \omega_0,Q(\omega_0,\mu_0) \bigr)$ is $(\tau,\gamma)$-Dio\-phan\-tine. The system~\eqref{eqfin} without the terms $f$, $g$, $h$ (the unperturbed system) at the parameter values $\omega=\omega_0$, $\sigma=0$, $\mu=\mu_0$ admits the reducible invariant $n$-torus $\{y=0, \, z=0\}$ with frequency vector $\omega_0$ and $(m+2p)\times(m+2p)$ Floquet matrix
\[
\begin{pmatrix} 0 & 0 \\ 0 & Q(\omega_0,\mu_0) \end{pmatrix}.
\]
According to Theorem~\ref{main}, the perturbed system~\eqref{eqfin} at the \emph{shifted} parameter values~\eqref{shift} has the reducible invariant $n$-torus $\{\by=0, \, \bz=0\}$ with \emph{the same} frequency vector and Floquet matrix. This torus is analytic (since the mappings~\eqref{mapp} are analytic in $\bx$) and depends on $\omega_0$ and $\mu_0$ in a $C^\infty$ way in the sense of Whitney (because the mappings~\eqref{mapp} are $C^\infty$ in $\omega$ and $\mu$). Thus, all the perturbed tori constitute a Whitney $C^\infty$-smooth family. This is a typical situation for the parametric KAM theory \cite{BHT90,BHS96Gro,BHS96LNM,BHN07,BH95,BCHV09}. The unperturbed and perturbed $n$-tori in the present setting pertain to the reversible context~2 since their phase space codimension is equal to $m+2p$, $\dim\Fix G=p$, and $m\geq 1$. The analogue of Theorem~\ref{main} for the reversible context~1 is Theorem~1.8 in \cite{BHN07}.

\begin{remark}
The $G$-reversibility conditions for the systems~\eqref{eqfin} are as follows (omitting the arguments $\omega$, $\sigma$, $\mu$): $QR=-RQ$ and
\begin{align*}
\xi(-y,Rz) &\equiv \xi(y,z), &\quad f(-x,-y,Rz) &\equiv f(x,y,z), \\
\eta(-y,Rz) &\equiv \eta(y,z), &\quad g(-x,-y,Rz) &\equiv g(x,y,z), \\
\zeta(-y,Rz) &\equiv -R\zeta(y,z), &\quad h(-x,-y,Rz) &\equiv -Rh(x,y,z).
\end{align*}
The condition that the changes of variables~\eqref{change} commute with $G$ (i.e., involution $G$ in the new coordinates has the original form $G:(\bx,\by,\bz)\mapsto(-\bx,-\by,R\bz)$) is equivalent to the identities (omitting the arguments $\omega$, $\mu$)
\begin{align*}
& & a(-\bx) &\equiv -a(\bx), \\
b^0(-\bx) &\equiv -b^0(\bx), &\quad b^1(-\bx) &\equiv b^1(\bx), &\quad b^2(-\bx)R &\equiv -b^2(\bx), \\
c^0(-\bx) &\equiv Rc^0(\bx), &\quad c^1(-\bx) &\equiv -Rc^1(\bx), &\quad c^2(-\bx)R &\equiv Rc^2(\bx).
\end{align*}
\end{remark}

\begin{remark}
Of course, the neighborhood $\cC$ in Theorem~\ref{main} contains the set
\begin{equation}
\bigl\{ (x,0,0,\omega,0,\mu) \bigm| x\in\mT^n, \; (\omega,\mu)\in\Gamma \bigr\}.
\label{vnutri}
\end{equation}
In fact, Theorem~\ref{main} is probably the simplest version of the statement on quasi-periodic stability of invariant tori in the reversible context~2. The first sentence of Theorem~\ref{main} can be sharpened as follows: ``There exists a closed set $\Gamma\subset\mR^{n+s}$ that is diffeomorphic to an $(n+s)$-di\-men\-sional ball, contains the point $(\omega_\ast,0)$ in its interior, and is such that for any complex neighborhood $\cC \subset (\mC/2\pi\mZ)^n\times\mC^{n+2m+2p+s}$ of the set~\eqref{vnutri}, there is a number $\rho>0$ with the following property.'' A similar remark refers to Theorem~\ref{source} in the next section.
\end{remark}

\section{The BCHV Theorem}
\label{sourcesec}

As was pointed out in Section~\ref{subject}, our proof of Theorem~\ref{main} is based on the results of the paper \cite{BCHV09}. Broer et al.\ \cite{BCHV09} deal with systems of the form~\eqref{eq1} that are reversible with respect to involution $G$~\eqref{rev1} and furthermore equivariant with respect to some action of the cyclic group $\mZ_l$ of $l\geq 1$ elements. We will need a particular case of the BCHV theorem where $l=1$ and the ``action-like'' variable $Y$ in Eq.~\eqref{eq1} is absent ($m=0$).

Consider an $(n+s)$-parameter analytic family of analytic systems
\begin{equation}
\begin{aligned}
\dot{x} &= \omega+\xi(z,\omega,\mu)+f(x,z,\omega,\mu), \\
\dot{z} &= Q(\omega,\mu)z+\zeta(z,\omega,\mu)+h(x,z,\omega,\mu)
\end{aligned}
\label{eqthey}
\end{equation}
reversible with respect to the involution $G:(x,z)\mapsto(-x,Rz)$. Here $x\in\mT^n$ and $z\in\cO_{2p}(0)$ are the phase space variables, $\omega\in\cO_n(\omega_\ast)$ and $\mu\in\cO_s(0)$ are external parameters ($n$, $p$, $s$ being non-negative integers), $R\in\GL(2p,\mR)$ is an involutive matrix with eigenvalues $1$ and $-1$ of multiplicity $p$ each, $Q$ is a $2p\times 2p$ matrix-valued function, $RQ(\omega,\mu)\equiv-Q(\omega,\mu)R$, $\xi=O(z)$, and $\zeta=O_2(z)$, whereas $f$ and $h$ are small perturbation terms.

The matrix-valued function $Q$ is assumed to satisfy the following two conditions. First, $\ker Q(\omega_\ast,0)\subset\Fix(-R)$, where $\Fix R$ and $\Fix(-R)$ are the $1$-eigen\-space and the $(-1)$-eigen\-space of the linear involution $R$, respectively. Second, the mapping $\mu\mapsto Q(\omega_\ast,\mu)\in\gl_{-R}$ is a \emph{versal unfolding} of the matrix $Q(\omega_\ast,0)$ with respect to the adjoint action of the group $\GL_{+R}$. Under these conditions, the following statement holds \cite{BCHV09}.

\begin{thrm}\label{source}
There exist a closed set $\Gamma\subset\mR^{n+s}$ that is diffeomorphic to an $(n+s)$-di\-men\-sional ball and contains the point $(\omega_\ast,0)$ in its interior, a number $\rho>0$, and a complex neighborhood $\cC \subset (\mC/2\pi\mZ)^n\times\mC^{n+2p+s}$ of the set
\[
\mT^n\times\{0\in\mR^{2p}\} \times \{\omega_\ast\}\times\{0\in\mR^s\}
\]
with the following property. For any $\tau>n-1$, $\gamma>0$, $L\in\mN$, and $\vare>0$ there is $\delta>0$ such that the following is valid. Suppose that the perturbation terms $f$ and $h$ can be holomorphically continued to the neighborhood $\cC$ and $|f|<\delta$, $|h|<\delta$ in $\cC$. Then there exist mappings
\begin{equation}
\begin{gathered}
a:\mT^n\times\Gamma\to\mR^n, \quad b:\mT^n\times\Gamma\to\mR^{2p}, \quad c:\mT^n\times\Gamma\to\gl(2p,\mR), \\
u:\Gamma\to\mR^n, \quad w:\Gamma\to\mR^s
\end{gathered}
\label{mapthey}
\end{equation}
possessing the following properties. First, the mappings~\eqref{mapthey} are analytic in $\bx\in\mT^n$ and $C^\infty$-smooth in $(\omega,\mu)\in\Gamma$. Second, all the partial derivatives of each component of these mappings of any order from $0$ to $L$ are everywhere smaller than $\vare$ in absolute value. Third, for any $(\omega,\mu)\in\Gamma$, each component of $a$, $b$, and $c$ as a function of $\bx$ can be holomorphically continued to the strip~\eqref{strip} and is smaller than $\vare$ in absolute value in this strip. Fourth, for any point $(\omega_0,\mu_0)\in\Gamma$, the nearly identical change of variables
\begin{equation}
x = \bx+a(\bx,\omega_0,\mu_0), \quad z = \bz+b(\bx,\omega_0,\mu_0)+c(\bx,\omega_0,\mu_0)\bz
\label{changthey}
\end{equation}
with $\bx\in\mT^n$ and $\bz\in\cO_{2p}(0)$ commutes with $G$. Fifth, for any point $(\omega_0,\mu_0)\in\Gamma$ such that the pair $\bigl( \omega_0,Q(\omega_0,\mu_0) \bigr)$ is $(\tau,\gamma)$-Dio\-phan\-tine, the system~\eqref{eqthey} at the parameter values
\[
\omega=\omega_0+u(\omega_0,\mu_0), \quad \mu=\mu_0+w(\omega_0,\mu_0)
\]
after the coordinate change~\eqref{changthey} takes the form
\[
\dot{\bx}=\omega_0+O(\bz), \quad \dot{\bz}=Q(\omega_0,\mu_0)\bz+O_2(\bz).
\]
\end{thrm}

The main novelty of the BCHV theorem compared with Theorem~1.8 in \cite{BHN07} is that the matrices $Q$ are allowed to be singular: the nondegeneracy condition $\det Q\neq 0$ is replaced with the much weaker condition $\ker Q\subset\Fix(-R)$. Since the operators $Q$ and $R$ anti-commute, one has
\[
Q(\Fix R)\subset\Fix(-R), \qquad Q\bigl(\Fix(-R)\bigr)\subset\Fix R,
\]
and $\ker Q$ is invariant under $R$, so that
\[
\ker Q = \bigl(\ker Q\cap\Fix R\bigr) \oplus \bigl(\ker Q\cap\Fix(-R)\bigr).
\]
Consequently, the condition $\ker Q\subset\Fix(-R)$ is equivalent to that $\ker Q\cap\Fix R = \{0\}$. This condition is not standard in KAM theory, and we will demonstrate its role with two ``toy'' examples.

The first example concerns the persistence of equilibria on the plane where $n=0$ and $p=s=1$. Consider a family of systems
\begin{equation}
\dot{z}_1 = z_2+\psi_1(z_1^2,z_2), \qquad \dot{z}_2 = \mu z_1+z_1\psi_2(z_1^2,z_2)
\label{ex1}
\end{equation}
reversible with respect to the involution $R:(z_1,z_2)\mapsto(-z_1,z_2)$. Here $z=(z_1,z_2)\in\cO_2(0)$ is the phase space variable, $\mu\in\cO_1(0)$ is a parameter, $\psi_1$, $\psi_2$ are arbitrary small functions, and $Q(\mu)=\left(\begin{smallmatrix} 0 & 1 \\ \mu & 0 \end{smallmatrix}\right)$ is a miniversal unfolding of $Q(0)=\left(\begin{smallmatrix} 0 & 1 \\ 0 & 0 \end{smallmatrix}\right)$ with respect to the adjoint action of the group $\GL_{+R}$ \cite{BCHV09,S92,S93,I96}. It is obvious that
\[
\ker Q(0) = \bigl\{ (z_1,0) \bigm| z_1\in\mR \bigr\} = \Fix(-R).
\]
We are looking for an equilibrium $(0,\fz)\in\Fix R$ of Eq.~\eqref{ex1} at $\mu=\fw$ where the linearization is similar to $Q(0)$. A point $(0,\fz)$ is an equilibrium of Eq.~\eqref{ex1} if and only if $\fz+\psi_1(0,\fz)=0$. This equation determines $\fz$ close to $0$ uniquely. The linearization matrix of Eq.~\eqref{ex1} around the equilibrium $(0,\fz)$ we have found is equal to
\[
\begin{pmatrix}
0 & 1+\partial\psi_1(0,\fz)/\partial z_2 \\ \mu+\psi_2(0,\fz) & 0
\end{pmatrix}.
\]
This matrix is similar to $Q(0)$ if and only if $\mu=\fw=-\psi_2(0,\fz)$. It is easy to verify that at this value of the parameter $\mu$, the system~\eqref{ex1} after the coordinate change
\[
z_1=\bz_1+\frac{\partial\psi_1(0,\fz)}{\partial z_2}\bz_1, \qquad z_2=\bz_2+\fz
\]
(which commutes with $R$) takes the form $\dot{\bz}_1=\bz_2+O_2(\bz_1,\bz_2)$, $\dot{\bz}_2=O_2(\bz_1,\bz_2)$.

Now consider a family of systems
\begin{equation}
\dot{z}_1 = z_2+z_2\psi_1(z_1,z_2^2), \qquad \dot{z}_2 = \mu z_1+\psi_2(z_1,z_2^2)
\label{ex2}
\end{equation}
reversible with respect to the involution $R:(z_1,z_2)\mapsto(z_1,-z_2)$. Here the matrices $Q(\mu)$ are the same as in Eq.~\eqref{ex1} and again constitute a miniversal unfolding of $Q(0)$ with respect to the adjoint action of $\GL_{+R}$, but $\ker Q(0)=\Fix R$. Again, we are looking for an equilibrium $(\fz,0)\in\Fix R$ of Eq.~\eqref{ex2} at $\mu=\fw$ where the linearization is similar to $Q(0)$. A point $(\fz,0)$ is an equilibrium of Eq.~\eqref{ex2} at $\mu=\fw$ if and only if $\fw\fz+\psi_2(\fz,0)=0$. The linearization matrix of Eq.~\eqref{ex2} around this equilibrium is equal to
\[
\begin{pmatrix}
0 & 1+\psi_1(\fz,0) \\ \fw+\partial\psi_2(\fz,0)/\partial z_1 & 0
\end{pmatrix}.
\]
Thus, for $\fz$ and $\fw$ we obtain the system of equations
\[
\fw\fz+\psi_2(\fz,0)=0, \qquad \fw+\partial\psi_2(\fz,0)/\partial z_1=0
\]
which has no solutions even in the simplest situation where $\psi_2$ is a non-zero constant.

Our second example is even more illuminative and concerns inhomogeneous linear systems. Consider a family of systems
\begin{equation}
\dot{z}=Q(\mu)z+\Psi(\mu)
\label{Psi}
\end{equation}
reversible with respect to a linear involution $R$, where $z\in\mR^N$, $\mu\in\cO_s(0)$, $Q$ is an $N\times N$ matrix-valued function, and $\Psi$ is an $N$-di\-men\-sional vector-valued function. The reversibility condition for Eq.~\eqref{Psi} is that $Q(\mu)\in\gl_{-R}$ and $\Psi(\mu)\in\Fix(-R)$ for each $\mu$. We are looking for a coordinate change $z=z'+\Delta(\mu)$ commuting with $R$ and reducing Eq.~\eqref{Psi} to the form $\dot{z}'=Q(\mu)z'$. Here $\Delta$ is also an $N$-di\-men\-sional vector-valued function, and a coordinate change $z=z'+\Delta(\mu)$ commutes with $R$ if and only if $\Delta(\mu)\in\Fix R$ for each $\mu$. For $z'$ we have the equation
\[
\dot{z}'=Q(\mu)\bigl( z'+\Delta(\mu) \bigr)+\Psi(\mu),
\]
whence $Q\Delta\equiv-\Psi$. Thus, a suitable $\Delta\in\Fix R$ can be found for any $\Psi\in\Fix(-R)$ if and only if the linear mapping $Q:\Fix R\to\Fix(-R)$ is an \emph{epimorphism}. If $N$ is even and $R$ has eigenvalues $1$ and $-1$ of multiplicity $N/2$ each (so that $\dim\Fix R=\dim\Fix(-R)=N/2$), this condition boils down to that $\ker Q\subset\Fix(-R)$.

If $Q(\mu)$ lies in the orbit $\fB\bigl( Q(0) \bigr)$ of $Q(0)$ under the adjoint action of the group $\GL_{+R}$ for each $\mu$, then one can reduce the equation $\dot{z}'=Q(\mu)z'$ to $\dot{\bz}=Q(0)\bz$ by an additional coordinate change $z'=A(\mu)\bz$, where $A:\cO_s(0)\to\GL_{+R}$ is a smooth $N\times N$ matrix-valued function and $A(0)$ is the $N\times N$ identity matrix. If $Q(\mu)$ is a versal unfolding of $Q(0)$ with respect to the adjoint action of $\GL_{+R}$, then for any small $N\times N$ matrix-valued function $\fQ:\cO_s(0)\to\gl_{-R}$ there exists a small value $\mu_\star\in\cO_s(0)$ of the parameter $\mu$ such that $Q(\mu_\star)+\fQ(\mu_\star) \in \fB\bigl( Q(0) \bigr)$.

\begin{remark}
Instead of the condition $\ker Q(\omega_\ast,0)\subset\Fix(-R)$, Broer et al.\ \cite{BCHV09} impose the following nondegeneracy condition on the systems~\eqref{eqthey}:
\[
\bigl[ \omega_\ast\partial/\partial x + Q(\omega_\ast,0)z\partial/\partial z, \: V \bigr]\neq 0
\]
whenever $V\in\cB^+\setminus\{0\}$, where $\cB^+$ is the space of ``constant'' $G$-equivariant vector fields and $[{\cdot},{\cdot}]$ is the Poisson bracket. The meaning of the word ``constant'' is not made precise explicitly in \cite{BCHV09} but it is clear from the text that one has in view vector fields of the form $\fb\partial/\partial z$ where $\fb\in\mR^{2p}$ is a constant vector (Prof.\ Han{\ss}mann has confirmed this in a private communication to me), see also \cite{H11RCD,H11IM}. Since
\[
[\omega\partial/\partial x + Qz\partial/\partial z, \: \fb\partial/\partial z] = -(Q\fb)\partial/\partial z,
\]
the condition in \cite{BCHV09} just mentioned means that $\ker Q(\omega_\ast,0)\cap\Fix R = \{0\}$, i.e., $\ker Q(\omega_\ast,0)\subset\Fix(-R)$. What is really needed in the proof of the BCHV theorem in \cite{BCHV09} is the equality $Q(\omega_\ast,0)(\Fix R)=\Fix(-R)$: the linear mapping $Q(\omega_\ast,0):\Fix R\to\Fix(-R)$ should be an epimorphism.
\end{remark}

If $\ker Q(\omega_\ast,0)\subset\Fix R$, one generically expects a reversible quasi-periodic center-saddle bifurcation in family~\eqref{eqthey} to occur \cite{H11RCD}. The Hamiltonian counterpart of this bifurcation scenario has been well studied \cite{H98,H07}.

\section{Proof of Theorem~\ref{main}}
\label{proof}

\subsection{The Crucial Trick}
\label{key}

Our goal is to reduce Theorem~\ref{main} to Theorem~\ref{source}. To this end we will treat the parameter $\sigma\in\cO_m(0)$ in the systems~\eqref{eqfin} as an additional phase space variable. In our paper \cite[Section~4.3]{S11}, such a strategy was called ``a naive approach to the reversible context~2'', and it was explained that this approach fails because systems~\eqref{eqfin} augmented by the equation $\dot{\sigma}=0$ are very strongly degenerate along the new ``normal'' variables $(y,\sigma,z)$. The key idea that enables one to overcome this difficulty is to replace the equation $\dot{\sigma}=0$ by the equation $\dot{\sigma}=\Lambda y$ where $\Lambda$ is a \emph{new additional external parameter} ranging in a neighborhood of the origin of the space $\gl(m,\mR)$ of $m\times m$ real matrices. To be more precise, instead of the $(n+m+s)$-parameter family~\eqref{eqfin} of systems with a phase space of dimension $n+m+2p$, we will consider the $(n+s+m^2)$-parameter family
\begin{equation}
\begin{aligned}
\dot{x} &= \omega+\xi(y,z,\omega,\sigma,\mu)+f(x,y,z,\omega,\sigma,\mu), \\
\dot{y} &= \sigma+\eta(y,z,\omega,\sigma,\mu)+g(x,y,z,\omega,\sigma,\mu), \\
\dot{\sigma} &= \Lambda y, \\
\dot{z} &= Q(\omega,\mu)z+\zeta(y,z,\omega,\sigma,\mu)+h(x,y,z,\omega,\sigma,\mu)
\end{aligned}
\label{eqnew}
\end{equation}
of systems with a phase space of dimension $n+2(m+p)$, where now the phase space variables are $(x,y,\sigma,z)$ and the external parameters are $(\omega,\mu,\Lambda)$. For $\Lambda=0$ we get the original family~\eqref{eqfin}. It is clear that each system in the family~\eqref{eqnew} is reversible with respect to the involution
\begin{equation}
\cG:(x,y,\sigma,z)\mapsto(-x,-y,\sigma,Rz).
\label{revnew}
\end{equation}
In particular, the $2(m+p)\times 2(m+p)$ matrices
\begin{equation}
\cQ(\omega,\mu,\Lambda) = \begin{pmatrix}
0 & I_m & 0 \\ \Lambda & 0 & 0 \\ 0 & 0 & Q(\omega,\mu)
\end{pmatrix}
\qquad \text{and} \qquad
\cR = \begin{pmatrix}
-I_m & 0 & 0 \\ 0 & I_m & 0 \\ 0 & 0 & R
\end{pmatrix}
\label{matrices}
\end{equation}
anti-commute for any $\omega$, $\mu$, and $\Lambda$ (here and henceforth, $I_m$ denotes the $m\times m$ identity matrix).

\subsection{Main Observations Concerning the Augmented Systems}
\label{observe}

First of all, note that since $R$ is an involutive matrix with eigenvalues $1$ and $-1$ of multiplicity $p$ each, the matrix $\cR$~\eqref{matrices} is involutive with eigenvalues $1$ and $-1$ of multiplicity $m+p$ each.

Second, since $\det Q(\omega_\ast,0)\neq 0$, we get
\begin{align*}
\ker\cQ(\omega_\ast,0,0) &= \bigl( \mR^m\times\{0\in\mR^m\}\times\{0\in\mR^{2p}\} \bigr) \\
{} &\subset \bigl( \mR^m\times\{0\in\mR^m\}\times\Fix(-R) \bigr) = \Fix(-\cR).
\end{align*}

Third, consider the anti-commuting $2m\times 2m$ matrices
\[
\cL(\Lambda) = \begin{pmatrix} 0 & I_m \\ \Lambda & 0 \end{pmatrix}
\qquad \text{and} \qquad
\cJ = \begin{pmatrix} -I_m & 0 \\ 0 & I_m \end{pmatrix},
\]
the matrix $\cJ=(-I_m)\oplus I_m$ being involutive with eigenvalues $1$ and $-1$ of multiplicity $m$ each. The group $\GL_{+\cJ}$ of non-singular $2m\times 2m$ real matrices commuting with $\cJ$ is $\bigl\{ \left(\begin{smallmatrix} A & 0 \\ 0 & B \end{smallmatrix}\right) \bigm| A,B\in\GL(m,\mR) \bigr\}$. Therefore, the orbit of $\cL(0)$ under the adjoint action of $\GL_{+\cJ}$ is
\[
\begin{pmatrix} 0 & \GL(m,\mR) \\ 0 & 0 \end{pmatrix}.
\]
Moreover, the space $\gl_{-\cJ}$ of $2m\times 2m$ real matrices anti-commuting with $\cJ$ is $\bigl\{ \left(\begin{smallmatrix} 0 & C \\ D & 0 \end{smallmatrix}\right) \bigm| C,D\in\gl(m,\mR) \bigr\}$. One concludes that the mapping $\Lambda\mapsto\cL(\Lambda)\in\gl_{-\cJ}$ is a miniversal unfolding of $\cL(0)$ with respect to the adjoint action of $\GL_{+\cJ}$. On the other hand, the mapping $\mu\mapsto Q(\omega_\ast,\mu)\in\gl_{-R}$ is a versal unfolding of $Q(\omega_\ast,0)$ with respect to the adjoint action of $\GL_{+R}$. The nilpotent matrix $\cL(0)$ and the non-singular matrix $Q(\omega_\ast,0)$ have no eigenvalues in common. Consequently \cite{S92,S93,I96}, the mapping
\[
(\mu,\Lambda) \mapsto \cL(\Lambda)\oplus Q(\omega_\ast,\mu) = \cQ(\omega_\ast,\mu,\Lambda) \in \gl_{-(\cJ\oplus R)} = \gl_{-\cR}
\]
is a versal unfolding of $\cQ(\omega_\ast,0,0)$ with respect to the adjoint action of $\GL_{+\cR}$.

\begin{remark}
We have found a miniversal unfolding of $\cL(0)$ with respect to the adjoint action of $\GL_{+\cJ}$ by extremely simple straightforward arguments. It is instructive to show how the same result follows from the general theorems of the papers \cite{S92,S93,I96} describing miniversal unfoldings of arbitrary infinitesimally reversible matrices. The Jordan normal form of $\cL(0)$ is the direct sum $\widetilde{\cL}(0)$ of $m$ nilpotent $2\times 2$ Jordan blocks $\left(\begin{smallmatrix} 0 & 1 \\ 0 & 0 \end{smallmatrix}\right)$. Let
\[
\widetilde{\cJ} = \diag(-1,1,-1,1,\ldots,-1,1) \in \GL(2m,\mR),
\]
then $\widetilde{\cL}(0)\widetilde{\cJ}=-\widetilde{\cJ}\widetilde{\cL}(0)$. According to Lemma~17 in \cite{S92} and Theorem~3(i) in \cite{S93} (see Figure~B in \cite{S93}), a miniversal unfolding of $\widetilde{\cL}(0)\in\gl_{-\widetilde{\cJ}}$ with respect to the adjoint action of $\GL_{+\widetilde{\cJ}}$ can be chosen to be the following matrix $\widetilde{\cL}(\lambda_{11},\ldots,\lambda_{mm})$. In the \mbox{$(2i-1)$-th} line of this matrix ($1\leq i\leq m$), the \mbox{$2i$-th} entry is equal to $1$ while all the other entries vanish. The \mbox{$2i$-th} line has the form
\[
(\lambda_{i1},0,\lambda_{i2},0,\ldots,\lambda_{im},0),
\]
where $\lambda_{ij}$ ($1\leq i\leq m$ and $1\leq j\leq m$) are independent parameters of the unfolding. After some analysis, one may verify that the paper \cite{I96} gives the same miniversal unfolding of $\widetilde{\cL}(0)$, see Corollary~3 and Table~VII in \cite{I96}. Now consider the $2m\times 2m$ matrix $\cS$ whose entries are all equal to zero except for the $2m$ entries
\[
\cS_{2i-1,i}=\cS_{2i,m+i}=1, \qquad 1\leq i\leq m.
\]
It is obvious that $\cS$ is non-singular (in fact, it is not hard to prove that $\det\cS=(-1)^{(m-1)m/2}$), and an easy calculation shows that
\[
\widetilde{\cJ}\cS=\cS\cJ, \qquad \widetilde{\cL}(\lambda_{11},\ldots,\lambda_{mm})\cS\equiv\cS\cL(\Lambda),
\]
where $\Lambda=(\lambda_{ij})_{1\leq i,j\leq m}\in\gl(m,\mR)$.
\end{remark}

We arrive at the conclusion that the family~\eqref{eqnew} of systems reversible with respect to the involution $\cG$~\eqref{revnew} satisfies all the conditions of Theorem~\ref{source} where
\begin{align*}
& \text{$m+p$ plays the role of $p$}, \\
& \text{$s+m^2$ plays the role of $s$}, \\
& \text{$(y,\sigma,z)$ plays the role of $z$}, \\
& \text{$(\mu,\Lambda)$ plays the role of $\mu$}, \\
& \text{$\cQ(\omega,\mu,\Lambda)$ plays the role of $Q(\omega,\mu)$}, \\
& \text{$(\eta,0,\zeta)$ with $0\in\mR^m$ plays the role of $\zeta$}, \\
& \text{$(g,0,h)$ with $0\in\mR^m$ plays the role of $h$}, \\
& \text{$\cR$ plays the role of $R$}, \\
& \text{$\cG$ plays the role of $G$}.
\end{align*}

\subsection{Consequences of the BCHV Theorem}
\label{apply}

Now we can apply Theorem~\ref{source} to the family~\eqref{eqnew}. The external parameters of this family are $(\omega,\mu,\Lambda)$, and Theorem~\ref{source} provides us with a closed set $\Gbig\subset\mR^{n+s}\times\gl(m,\mR)$ such that (i) $\Gbig$ is diffeomorphic to an $(n+s+m^2)$-di\-men\-sional ball; (ii) $\Gbig$ contains the point $(\omega_\ast,0,0)$ in its interior; (iii) there are dynamical consequences for all the points $(\omega_0,\mu_0,\Lambda_0)\in\Gbig$ such that the pair $\bigl( \omega_0,\cQ(\omega_0,\mu_0,\Lambda_0) \bigr)$ is $(\tau,\gamma)$-Dio\-phan\-tine. However, for our purposes it will be sufficient to confine ourselves with the case $\Lambda_0=0$. The spectrum of the matrix $\cQ(\omega,\mu,0)$ is just the spectrum of the matrix $Q(\omega,\mu)$ plus $2m$ zero eigenvalues. Therefore, the pair $\bigl( \omega_0,\cQ(\omega_0,\mu_0,0) \bigr)$ is $(\tau,\gamma)$-Dio\-phan\-tine if and only if the pair $\bigl( \omega_0,Q(\omega_0,\mu_0) \bigr)$ is $(\tau,\gamma)$-Dio\-phan\-tine.

Thus, according to Theorem~\ref{source}, there exist a closed set $\Gamma\subset\mR^{n+s}$ that is diffeomorphic to an $(n+s)$-di\-men\-sional ball and contains the point $(\omega_\ast,0)$ in its interior, a number $\rho>0$, and a complex neighborhood $\cC \subset (\mC/2\pi\mZ)^n\times\mC^{n+2m+2p+s}$ of the set~\eqref{chord} with the following property. For any $\tau>n-1$, $\gamma>0$, $L\in\mN$, and $\vare>0$ there is $\delta>0$ such that the following holds. Suppose that the perturbation terms $f$, $g$, $h$ can be holomorphically continued to the neighborhood $\cC$ and $|f|<\delta$, $|g|<\delta$, $|h|<\delta$ in $\cC$. Then there exist mappings
\begin{equation}
\begin{gathered}
a:\mT^n\times\Gamma\to\mR^n, \\
\begin{aligned}
b^0&:\mT^n\times\Gamma\to\mR^m, &\quad b^1,b^2&:\mT^n\times\Gamma\to\gl(m,\mR), &\quad b^3&:\mT^n\times\Gamma\to\mR^{m\times 2p}, \\
c^0&:\mT^n\times\Gamma\to\mR^m, &\quad c^1,c^2&:\mT^n\times\Gamma\to\gl(m,\mR), &\quad c^3&:\mT^n\times\Gamma\to\mR^{m\times 2p}, \\
d^0&:\mT^n\times\Gamma\to\mR^{2p}, &\quad d^1,d^2&:\mT^n\times\Gamma\to\mR^{2p\times m}, &\quad d^3&:\mT^n\times\Gamma\to\gl(2p,\mR),
\end{aligned} \\
u:\Gamma\to\mR^n, \quad v:\Gamma\to\mR^s, \quad W:\Gamma\to\gl(m,\mR)
\end{gathered}
\label{mapwe}
\end{equation}
possessing the following properties. First, the mappings~\eqref{mapwe} are analytic in $\bx\in\mT^n$ and $C^\infty$-smooth in $(\omega,\mu)\in\Gamma$. Second, all the partial derivatives of each component of these mappings of any order from $0$ to $L$ are everywhere smaller than $\vare$ in absolute value. Third, for any $(\omega,\mu)\in\Gamma$, each component of $a$, $b^j$, $c^j$, $d^j$ ($0\leq j\leq 3$) as a function of $\bx$ can be holomorphically continued to the strip~\eqref{strip} and is smaller than $\vare$ in absolute value in this strip. Fourth, for any point $(\omega_0,\mu_0)\in\Gamma$, the nearly identical change of variables
\begin{equation}
\begin{aligned}
x &= \bx+a(\bx,\omega_0,\mu_0), \\
y &= \by+b^0(\bx,\omega_0,\mu_0)+b^1(\bx,\omega_0,\mu_0)\by+b^2(\bx,\omega_0,\mu_0)\bsigma+b^3(\bx,\omega_0,\mu_0)\bz, \\
\sigma &= \bsigma+c^0(\bx,\omega_0,\mu_0)+c^1(\bx,\omega_0,\mu_0)\by+c^2(\bx,\omega_0,\mu_0)\bsigma+c^3(\bx,\omega_0,\mu_0)\bz, \\
z &= \bz+d^0(\bx,\omega_0,\mu_0)+d^1(\bx,\omega_0,\mu_0)\by+d^2(\bx,\omega_0,\mu_0)\bsigma+d^3(\bx,\omega_0,\mu_0)\bz
\end{aligned}
\label{changwe}
\end{equation}
with $\bx\in\mT^n$, $\by\in\cO_m(0)$, $\bsigma\in\cO_m(0)$, $\bz\in\cO_{2p}(0)$ commutes with $\cG$. Fifth, for any point $(\omega_0,\mu_0)\in\Gamma$ such that the pair $\bigl( \omega_0,Q(\omega_0,\mu_0) \bigr)$ is $(\tau,\gamma)$-Dio\-phan\-tine, the system~\eqref{eqnew} at the parameter values
\[
\omega=\omega_0+u(\omega_0,\mu_0), \quad \mu=\mu_0+v(\omega_0,\mu_0), \quad \Lambda=W(\omega_0,\mu_0)
\]
after the coordinate change~\eqref{changwe} takes the form
\begin{equation}
\begin{gathered}
\dot{\bx}=\omega_0+O(\by,\bsigma,\bz), \\
\dot{\by}=\bsigma+O_2(\by,\bsigma,\bz), \quad \dot{\bsigma}=O_2(\by,\bsigma,\bz), \quad \dot{\bz}=Q(\omega_0,\mu_0)\bz+O_2(\by,\bsigma,\bz).
\end{gathered}
\label{normal}
\end{equation}

We do not assume the neighborhood $\cC$ to lie in $(\mC/2\pi\mZ)^n\times\mC^{n+2m+2p+s+m^2}$ because the perturbation terms $f$, $g$, $h$ do not depend on $\Lambda$.

\subsection{Absence of a Shift along the Parameter $\Lambda$}
\label{Wzero}

Assuming $\tau>n-1$, $\gamma>0$, and $L\in\mN$ to be fixed, consider an arbitrary point $(\omega_0,\mu_0)\in\Gamma$ such that the pair $\bigl( \omega_0,Q(\omega_0,\mu_0) \bigr)$ is $(\tau,\gamma)$-Dio\-phan\-tine. We will drop the arguments $\omega_0$, $\mu_0$ of functions $a$, $b^j$, $c^j$, $d^j$ ($0\leq j\leq 3$), $u$, $v$, $W$, and $Q$. Our first (and principal) aim is to prove that $W=0$ for $\vare$ sufficiently small taking into account the very simple form $\dot{\sigma}=\Lambda y$ of the equation for $\dot{\sigma}$ in the systems~\eqref{eqnew}.

Suppose that up to $O_2(\by,\bsigma,\bz)$, the equation for $\dot{\bx}$ in the normalized systems~\eqref{normal} has the form
\begin{equation}
\dot{\bx} = \omega_0+\chi^1(\bx)\by+\chi^2(\bx)\bsigma+\chi^3(\bx)\bz + O_2(\by,\bsigma,\bz)
\label{normaladd}
\end{equation}
with analytic coefficients
\begin{equation}
\chi^1:\mT^n\to\mR^{n\times m}, \quad \chi^2:\mT^n\to\mR^{n\times m}, \quad \chi^3:\mT^n\to\mR^{n\times 2p}
\label{chi}
\end{equation}
(of course, these coefficients depend also on $\omega_0$ and $\mu_0$). On $\mT^n$, each component of the mappings~\eqref{chi} does not exceed some $\vare$-independent constant $E>0$ in absolute value. By virtue of Eqs.~\eqref{changwe}--\eqref{normaladd}, the equation $\dot{\sigma}=\Lambda y$ with $\Lambda=W$ takes the form
\begin{equation}
\begin{aligned}
\frac{\partial c^0}{\partial\bx}\bigl( \omega_0+\chi^1\by+\chi^2\bsigma+\chi^3\bz \bigr) &+
\left( \frac{\partial c^1}{\partial\bx}\omega_0 \right)\by + c^1\bsigma +
\left( \frac{\partial c^2}{\partial\bx}\omega_0 \right)\bsigma \\
{} &+ \left( \frac{\partial c^3}{\partial\bx}\omega_0 \right)\bz +c^3Q\bz +
O_2(\by,\bsigma,\bz) \\
{} &= W\bigl( \by+b^0+b^1\by+b^2\bsigma+b^3\bz \bigr).
\end{aligned}
\label{ravno}
\end{equation}

Now we will need the following standard and easy lemma ubiquitous in the problems concerning small divisors.

\begin{lemma}\label{rho}
Suppose that two holomorphic functions $F,\Phi:\cU_n(\rho)\to\mC$ (see Eq.~\eqref{strip}) with zero average satisfy the identity
\[
\frac{\partial\Phi}{\partial\bx}\omega_0 \equiv F
\]
in $\cU_n(\rho)$ where the vector $\omega_0\in\mR^n$ is $(\tau,\gamma)$-Dio\-phan\-tine. Then for any number $\rho'$ in the interval $0<\rho'<\rho$, there holds the estimate
\begin{equation}
\sup_{\bx\in\cU_n(\rho')}\bigl| \Phi(\bx) \bigr| \leq
\frac{\fC_{n,\tau}}{\gamma(\rho-\rho')^{n+\tau}}
\sup_{\bx\in\cU_n(\rho)}\bigl| F(\bx) \bigr|,
\label{estimate}
\end{equation}
where $\fC_{n,\tau}>0$ is a certain constant depending on $n$ and $\tau$ only.
\end{lemma}

\begin{remark}
In fact, a sharper result is valid, with $(\rho-\rho')^\tau$ in place of $(\rho-\rho')^{n+\tau}$ in Eq.~\eqref{estimate}, see \cite[Lemma~3.15]{dL01} and \cite{R75}. However, for almost all the studies in KAM theory in the analytic category and, in particular, for the present paper, the estimate~\eqref{estimate} is enough.
\end{remark}

Introduce the notation
\[
\|W\| = \max_{i,j=1}^m |W_{ij}|
\]
and fix an arbitrary number $\rho'$ in the interval $0<\rho'<\rho$. Equating the constant terms in the left-hand side and right-hand side of Eq.~\eqref{ravno}, we see that
\[
\frac{\partial c^0}{\partial\bx}\omega_0=Wb^0.
\]
This equality holds for all $\bx\in\mT^n$ and, consequently, for all $\bx\in\cU_n(\rho)$. In the strip $\cU_n(\rho)$, each component of $b^0$ is less than $\vare$ in absolute value. According to Lemma~\ref{rho}, in the strip $\cU_n(\rho')$, each component of $c^0-\langle c^0\rangle$ does not exceed
\[
\frac{\fC_{n,\tau}m\|W\|\vare}{\gamma(\rho-\rho')^{n+\tau}}
\]
in absolute value, where $\langle c^0\rangle$ is the average of $c^0$ (such a notation will be also used below). The Cauchy estimate implies that for real $\bx\in\mT^n$, each entry of the $m\times n$ Jacobi matrix $\partial c^0/\partial\bx$ does not exceed
\begin{equation}
\frac{m\fC_{n,\tau}\|W\|\vare}{\gamma(\rho-\rho')^{n+\tau}\rho'}
\label{Jacobi}
\end{equation}
in absolute value.

Now equate the terms linear in $\by$ in the left-hand side and right-hand side of Eq.~\eqref{ravno}:
\[
\frac{\partial c^0}{\partial\bx}\chi^1 + \frac{\partial c^1}{\partial\bx}\omega_0 = W+Wb^1
\]
and average over $\bx\in\mT^n$:
\[
\left\langle \frac{\partial c^0}{\partial\bx}\chi^1 \right\rangle = W+W\langle b^1\rangle.
\]
For $b^1$ sufficiently small (i.e., for $\vare$ sufficiently small), the $m\times m$ matrix $I_m+\langle b^1\rangle$ is non-singular, and each entry of its inverse is less than $2$ in absolute value (instead of $2$, one may use any constant greater than $1$). Recalling that each component of $\chi^1$ does not exceed $E$ in absolute value and employing the estimate~\eqref{Jacobi} for $\partial c^0/\partial\bx$, we arrive at the conclusion that the equality
\[
W = \left\langle \frac{\partial c^0}{\partial\bx}\chi^1 \right\rangle \bigl( I_m+\langle b^1\rangle \bigr)^{-1}
\]
implies that
\[
\|W\| \leq \frac{2Em^2n\fC_{n,\tau}\|W\|\vare}{\gamma(\rho-\rho')^{n+\tau}\rho'}.
\]
Consequently, $\|W\|=0$ for $\vare$ small enough.

\begin{remark}
One may wonder whether it is possible to deduce just from evenness arguments that $\bigl\langle (\partial c^0/\partial\bx)\chi^1 \bigr\rangle=0$. However, this is not the case. Indeed, $c^0(\bx)$ is even in $\bx$ because the coordinate change~\eqref{changwe} commutes with $\cG$ (see Eq.~\eqref{revnew}), and therefore $\partial c^0/\partial\bx$ is odd in $\bx$. On the other hand, $\chi^1(\bx)$ is also odd in $\bx$ because the systems~\eqref{normal} are reversible with respect to $\cG$.
\end{remark}

\begin{remark}
Another way to prove that $W$ vanishes is based on an additional coordinate change $\bx=\varp+\vart(\varp)\by$ where $\varp\in\mT^n$. Since $\omega_0$ is $(\tau,\gamma)$-Dio\-phan\-tine, one can remove the term linear in $\by$ in the equation for $\dot{\varp}$ by a suitable choice of the function $\vart:\mT^n\to\mR^{n\times m}$.
\end{remark}

\subsection{Completion of the Proof}
\label{complete}

Since $W=0$, the left-hand side of Eq.~\eqref{ravno} vanishes and its constant term $(\partial c^0/\partial\bx)\omega_0$ is zero. This implies that $c^0$ is independent of $\bx$ because the vector $\omega_0$ is non-resonant. Now one can consider the terms linear in $\by$, $\bsigma$, and $\bz$ in the left-hand side of Eq.~\eqref{ravno} and obtain that
\begin{gather}
(\partial c^1/\partial\bx)\omega_0 = 0,
\label{q1} \\
c^1+(\partial c^2/\partial\bx)\omega_0 = 0,
\label{q2} \\
(\partial c^3/\partial\bx)\omega_0+c^3Q = 0.
\label{q3}
\end{gather}
Eq.~\eqref{q1} implies that $c^1$ is independent of $\bx$. On the other hand, $\langle c^1\rangle=0$ according to Eq.~\eqref{q2}. Therefore, $c^1=0$, and Eq.~\eqref{q2} implies that $c^2$ is independent of $\bx$. Since the pair $(\omega_0,Q)$ is $(\tau,\gamma)$-Dio\-phan\-tine (in fact, non-resonance would be enough here), one can conclude from Eq.~\eqref{q3} that all the Fourier coefficients $c^3_k$ of $c^3$ with $k\in\mZ^n\setminus\{0\}$ vanish (i.e., $c^3$ is independent of $\bx$ as well). Moreover, Eq.~\eqref{q3} implies that $c^3_0Q=0$ whence $c^3_0=0$ (because $\det Q\neq 0$). Consequently, $c^3=0$.

Thus, $\sigma=c^0+(I_m+c^2)\bsigma$ where $c^0$ and $c^2$ are small and independent of $\bx$. Since $W=0$ and $\dot{\sigma}\equiv 0$, it follows that $\dot{\bsigma}\equiv 0$. The transformation~\eqref{changwe} casts the invariant plane $\{\sigma=c^0\}$ of the system~\eqref{eqnew} at the parameter values
\[
\omega=\omega_0+u(\omega_0,\mu_0), \quad \mu=\mu_0+v(\omega_0,\mu_0), \quad \Lambda=0
\]
to the invariant plane $\{\bsigma=0\}$ of the system~\eqref{normal}. The restriction of Eq.~\eqref{normal} to $\{\bsigma=0\}$ has the form~\eqref{goal}.

We have verified the equalities $W=0$, $\partial c^0/\partial\bx=0$, $c^1=0$, $\partial c^2/\partial\bx=0$, and $c^3=0$ (provided that $\vare$ is sufficiently small) for points $(\omega_0,\mu_0)\in\Gamma$ such that the pair $\bigl( \omega_0,Q(\omega_0,\mu_0) \bigr)$ is $(\tau,\gamma)$-Dio\-phan\-tine. But we are interested in the coordinate change~\eqref{changwe} for such points $(\omega_0,\mu_0)$ only. Therefore, one may set $W$, $\partial c^0/\partial\bx$, $c^1$, $\partial c^2/\partial\bx$, and $c^3$ to vanish everywhere. Note also that the coefficients $b^2$, $c^2$, and $d^2$ in Eq.~\eqref{changwe} are irrelevant as far as the plane $\{\bsigma=0\}$ is concerned.

We have arrived at the conclusion of Theorem~\ref{main}, with
\begin{align*}
& \text{$b^3$ playing the role of $b^2$}, \\
& \text{$d^0$, $d^1$, $d^3$ playing the roles of $c^0$, $c^1$, $c^2$, respectively}, \\
& \text{$c^0$ playing the role of $v$}, \\
& \text{$v$ playing the role of $w$}.
\end{align*}
The proof of Theorem~\ref{main} is completed.

\section{R\"ussmann Nondegeneracy}
\label{Ruesssec}

Most probably, Theorem~\ref{main} enables one to carry over to the reversible context~2 such phenomena as the so-called excitation of elliptic normal modes and partial preservation of the frequencies and Floquet exponents of the unperturbed tori. The excitation of elliptic normal modes is well known in the reversible context~1 (see \cite{BHS96LNM,S95Cha,S98} and references therein) and the Hamiltonian context (see \cite{AKN06,BS10,BHS96LNM} and references therein). The partial preservation of frequencies and Floquet exponents has been also studied in detail in the reversible context~1 and other ``conventional'' contexts of KAM theory \cite{S06,S07Stek}. We hope to explore these topics in subsequent publications. In the present paper, we only consider the simplest application of Theorem~\ref{main}, namely, an analogue of the R\"ussmann nondegeneracy condition for the reversible context~2 in the absence of the ``normal'' variable $z$. Our exposition will be less formal than that in Sections~\ref{mainsec}--\ref{proof}.

Consider a family of systems
\begin{equation}
\dot{x}=H(Y,\lambda), \qquad \dot{Y}=0
\label{simpleq}
\end{equation}
reversible with respect to the involution $G:(x,Y)\mapsto(-x,Y)$, where $x\in\mT^n$ and $Y\in\mR^m$ are the phase space variables and $\lambda\in\mR^s$ is an external parameter ($n\geq 1$, $m\geq 0$, $s\geq 0$, $m+s\geq 1$), cf.\ Eq.~\eqref{eq1}. We suppose that $(Y,\lambda)$ ranges in the closure $\fW\subset\mR^{m+s}$ of a bounded connected open domain in $\mR^{m+s}$ and the function $H:\fW\to\mR^n$ is analytic.

The family~\eqref{simpleq} of integrable dynamical systems is said to be \emph{KAM-stable} if any sufficiently small analytic $G$-reversible perturbation of Eq.~\eqref{simpleq} admits a Whitney smooth family of analytic invariant $n$-tori carrying quasi-periodic motions and close to the unperturbed tori $\{Y=\const\}$, the Lebesgue measure of the complement of the union of the perturbed tori in $\mT^n\times\fW$ tending to zero as the perturbation magnitude tends to $0$. More formal definitions of KAM-stability are presented in \cite{S95JDCS,S98}. It turns out that the family~\eqref{simpleq} is KAM-stable if and only if it is nondegenerate in the sense of R\"ussmann: the image $H(\fW)$ of the frequency mapping $H:\fW\to\mR^n$ does not lie in any linear hyperplane passing through the origin \cite{S95JDCS,BHS96LNM,S98,WX09} (here the analyticity of $H$ and the connectedness of $\fW$ are essential). The Hamiltonian counterpart of this theorem is very well known, see \cite{AKN06,BS10,S95JDCS,BHS96LNM,S98,H11IM,R01,R05} and references therein. One of the main ingredients of the proofs is the following number-theoretical lemma, see \cite{S95JDCS,BHS96LNM,S95Cha,R01} and references therein.

\begin{lemma}\label{plane}
Let $\fW\subset\mR^N$ be the closure of a bounded connected open domain in $\mR^N$ and let $H:\fW\to\mR^n$ be an analytic function. Assume that the image $H(\fW)$ of $H$ does not lie in any linear hyperplane passing through the origin. Then there exists a positive integer $r$ such that the following holds. For any fixed value $\tau>nr-1$, the Lebesgue measure of the set of points $a\in\fW$ for which the vector $H^\sharp(a)$ is \emph{not} $(\tau,\gamma)$-Dio\-phan\-tine tends to zero as $\gamma\to 0$ uniformly with respect to all the $C^r$-functions $H^\sharp:\fW\to\mR^n$ in some $C^r$-neighborhood of the function $H$.
\end{lemma}

The R\"ussmann nondegeneracy condition is very weak: for any dimension $n$ of the frequency space, it is easy to construct an analytic mapping $H:\cO_1(0)\to\mR^n$ whose image (a curve in $\mR^n$) does not lie in any linear hyperplane passing through the origin \cite{S95JDCS,BHS96LNM}. The simplest example is $H(Y)=(1,Y,Y^2,\ldots,Y^{n-1})$.

Now consider an $(m+s)$-parameter analytic family of analytic systems
\begin{equation}
\dot{x}=F(\sigma,\mu)+\xi(y,\sigma,\mu), \qquad \dot{y}=\sigma+\eta(y,\sigma,\mu)
\label{moreq}
\end{equation}
reversible with respect to the involution $G:(x,y)\mapsto(-x,-y)$, where $x\in\mT^n$, $y\in\cO_m(0)$ are the phase space variables, $\sigma\in\cO_m(0)$, $\mu\in\fW\subset\mR^s$ are external parameters, and $\xi=O_2(y)$, $\eta=O_2(y)$ ($n$, $m$, $s$ being positive integers and $\fW$ being the closure of a bounded connected open domain in $\mR^s$), cf.\ Eqs.~\eqref{eq2} and~\eqref{eqfin}. The $G$-reversibility condition consists in that the functions $\xi$ and $\eta$ are even in $y$. The systems~\eqref{eqfin} may suggest that the requirement $\xi=O(y)$ would be more natural than $\xi=O_2(y)$. However, these two requirements are equivalent in the case of Eq.~\eqref{moreq}: if $\xi=O(y)$ and $\xi$ is even in $y$ then $\xi=O_2(y)$. For $\sigma=0$, the $n$-torus $\{y=0\}$ is invariant under both the system~\eqref{moreq} and the involution $G$ and carries conditionally periodic motions with frequency vector $F(0,\mu)$.

\begin{defn}
The family~\eqref{moreq} pertaining to the reversible context~2 is said to be \emph{nondegenerate in the sense of R\"ussmann} if the set
\begin{equation}
\bigl\{ F(0,\mu) \bigm| \mu\in\fW \bigr\} \subset \mR^n
\label{curled}
\end{equation}
does not lie in any linear hyperplane passing through the origin.
\end{defn}

\begin{thrm}\label{Ruess}
Let the family~\eqref{moreq} be nondegenerate in the sense of R\"ussmann. Then for any sufficiently small analytic $G$-reversible perturbation
\begin{equation}
\begin{aligned}
\dot{x} &= F(\sigma,\mu)+\xi(y,\sigma,\mu)+f(x,y,\sigma,\mu), \\
\dot{y} &= \sigma+\eta(y,\sigma,\mu)+g(x,y,\sigma,\mu)
\end{aligned}
\label{perturb}
\end{equation}
of Eq.~\eqref{moreq}, there exist a subset $\fW_0\subset\fW$ and a small $C^\infty$-function $\Theta:\fW\to\mR^m$ such that the following holds. For any $\mu\in\fW_0$, the system~\eqref{perturb} with $\sigma=\Theta(\mu)$ admits an analytic invariant $n$-torus that is close to the torus $\{y=0\}$ and carries quasi-periodic motions with a Diophantine frequency vector close to $F(0,\mu)$. All such tori constitute a Whitney $C^\infty$-smooth family. The Lebesgue measure of $\fW\setminus\fW_0$ tends to zero as the perturbation magnitude tends to $0$.
\end{thrm}

Let us prove this theorem omitting some boring technical details. According to Lemma~\ref{plane}, there exists $r\in\mN$ possessing the following property. For any fixed value $\tau>nr-1$, the Lebesgue measure of the set of points $\mu\in\fW$ for which the vector $F^\sharp(\mu)$ is not $(\tau,\gamma)$-Dio\-phan\-tine tends to zero as $\gamma\to 0$ uniformly with respect to all the $C^r$-functions $F^\sharp:\fW\to\mR^n$ in some $C^r$-neighborhood of the mapping $\mu\mapsto F(0,\mu)$. For our purposes, it will suffice to deal with $C^\infty$-functions $F^\sharp$. Fix a number $\tau>nr-1$ and consider the family of $G$-reversible systems
\begin{equation}
\begin{aligned}
\dot{x} &= \omega+\xi(y,\sigma,\mu)+f(x,y,\sigma,\mu), \\
\dot{y} &= \sigma+\eta(y,\sigma,\mu)+g(x,y,\sigma,\mu)
\end{aligned}
\label{extend}
\end{equation}
instead of Eq.~\eqref{perturb}, where $\omega\in\mR^n$ is a \emph{new additional external parameter} ranging in some closed neighborhood $\fO$ of the set~\eqref{curled}.

For any $\gamma>0$, Theorem~\ref{main} provides us with small $C^\infty$-mappings
\[
u:\fO\times\fW\to\mR^n, \quad v:\fO\times\fW\to\mR^m, \quad w:\fO\times\fW\to\mR^s
\]
such that for any $(\tau,\gamma)$-Dio\-phan\-tine vector $\omega_0\in\fO$ and any point $\mu_0\in\fW$, the system~\eqref{extend} at the parameter values~\eqref{shift} has an analytic invariant $n$-torus (close to the torus $\{y=0\}$) carrying quasi-periodic motions with frequency vector $\omega_0$. Here we suppose that the perturbation terms $f$ and $g$ are small enough (the corresponding smallness requirement depends on $\gamma$) and $\mu_0+w(\omega_0,\mu_0)$ still lies in $\fW$.

The equation
\begin{equation}
F\bigl( v(\omega,\mu), \, \mu+w(\omega,\mu) \bigr) = \omega+u(\omega,\mu)
\label{Phi}
\end{equation}
(with $\mu\in\fW$) can be solved with respect to $\omega$:
\[
\omega=\Phi(\mu),
\]
where $\Phi:\fW\to\mR^n$ is a $C^\infty$-function close to the mapping $\mu\mapsto F(0,\mu)$. For any point $\mu_0\in\fW$ such that $\Phi(\mu_0)$ is $(\tau,\gamma)$-Dio\-phan\-tine, the \emph{original} perturbed system~\eqref{perturb} at the parameter values
\begin{equation}
\sigma=v\bigl(\Phi(\mu_0),\mu_0 \bigr), \qquad \mu=\mu_0+w\bigl(\Phi(\mu_0),\mu_0 \bigr)
\label{newshift}
\end{equation}
has an invariant $n$-torus carrying quasi-periodic motions with frequency vector $\Phi(\mu_0)$. Indeed, the system~\eqref{extend} at the parameter values~\eqref{shift} with $\omega_0=\Phi(\mu_0)$ \emph{coincides} with the system~\eqref{perturb} at the parameter values~\eqref{newshift}.

The equation
\begin{equation}
\mu=\mu_0+w\bigl(\Phi(\mu_0),\mu_0 \bigr)
\label{Upsilon}
\end{equation}
(with $\mu\in\fW$) can be solved with respect to $\mu_0$:
\[
\mu_0=\Upsilon(\mu),
\]
where $\Upsilon:\fW\to\mR^s$ is a $C^\infty$-function close to the identity mapping $\mu\mapsto\mu$. We have arrived at the following conclusion. For any point $\mu\in\fW$ such that
\[
F^\sharp(\mu) = \Phi\bigl( \Upsilon(\mu) \bigr) \in \mR^n
\]
is $(\tau,\gamma)$-Dio\-phan\-tine, the system~\eqref{perturb} for
\[
\sigma = \Theta(\mu) = v\bigl( F^\sharp(\mu),\Upsilon(\mu) \bigr)
\]
has an invariant $n$-torus carrying quasi-periodic motions with frequency vector $F^\sharp(\mu)$. The $C^\infty$-function $F^\sharp:\fW\to\mR^n$ is close to the mapping $\mu\mapsto F(0,\mu)$. The $C^\infty$-function $\Theta:\fW\to\mR^m$ is small. Let
\[
\fW_0 = \bigl\{ \mu\in\fW \bigm| \text{$F^\sharp(\mu)$ is $(\tau,\gamma)$-Diophantine} \bigr\},
\]
then the Lebesgue measure of $\fW\setminus\fW_0$ tends to zero as $\gamma\to 0$. Theorem~\ref{Ruess} has been proven. To deduce Theorem~\ref{Ruess} from Theorem~\ref{main} and Lemma~\ref{plane}, we just introduced an additional external parameter $\omega$ and invoked the implicit function theorem twice to solve Eqs.~\eqref{Phi} and~\eqref{Upsilon}.

The author does not know whether the nondegeneracy of the family~\eqref{moreq} in the sense of R\"ussmann is necessary for the conclusion of Theorem~\ref{Ruess}.

\begin{remark}
Introducing a new external parameter to achieve the ``maximal'' nondegeneracy of the systems was the key step in the proofs of both Theorem~\ref{main} and Theorem~\ref{Ruess}. This parameter was $\Lambda\in\gl(m,\mR)$ in the case of Theorem~\ref{main} (see Section~\ref{key}) and $\omega\in\mR^n$ in the case of Theorem~\ref{Ruess}. In KAM theory, such a trick with an additional parameter is usually called the \emph{Herman method} \cite{BS10,S95JDCS,BHS96Gro,BHS96LNM,S06,S07Stek,S07DCDS,S95Cha,KP01}; it was first presented by M.~R.~Herman in his talk at the International Conference on Dynamical Systems in Lyons in 1990 (see also \cite{Y92}). However, in the proofs of Theorem~\ref{main} and Theorem~\ref{Ruess}, we employed the additional external parameter in different ways. In the case of Theorem~\ref{main}, the principal part of the proof was to verify that the shift along $\Lambda$ vanishes for $\Lambda_0=0$ (see Section~\ref{Wzero}). In the case of Theorem~\ref{Ruess}, we just eliminated $\omega$ by solving Eq.~\eqref{Phi}. The latter course is typical for the Herman method.

By the way, there are many situations in KAM theory where various statements can be deduced in a very straightforward manner from statements pertaining to other setups. The Herman method is just one of the tools for such a reduction. In \cite{S12IM}, we presented a list of such situations which was of course not exhaustive. The settings missed in the list of \cite{S12IM} are exemplified by a reduction of the perturbative KAM theorems to the so-called a~posteriori ones, see e.g.\ \cite{HCLMF16,FdLS09}.
\end{remark}

\section{Concluding Remarks}
\label{conclu}

The ``splitting'' of the reversible KAM theory into the reversible contexts~1 and~2 is one of the phenomena indicating that the striking similarity between the Hamiltonian context and reversible context \cite{BHS96LNM} does not go very far. A detailed discussion of the Hamiltonian-reversible parallelism is beyond the present study (see \cite{S98} for some remarks), but we would like to mention one more example of important differences between the properties of quasi-periodic motions in Hamiltonian and reversible systems. Namely, a reversible analogue of the Nekhoroshev theory on the effective (exponential) stability of the action variables in analytic nearly integrable Hamiltonian systems \cite{N77,N79} is hardly possible, see \cite[\S~4.2.5]{BHS96LNM} and \cite[Section~2]{S98}. The reason is that within the reversible realm, there is no coupling between the resonant unperturbed frequency vectors and the averaged evolution of the action variables. Of course, this does not exclude effective stability in analytic nearly integrable reversible systems satisfying special conditions, e.g.\ in the case where the unperturbed system is a collection of harmonic oscillators with a Diophantine vector of frequencies \cite{DL99,Z05}.

The books \cite{AKN06,D14,BHS96LNM,A88,A04} contain brief accounts of the Nekhoroshev theory (see also the paper \cite{P00}); for a modern survey see e.g.\ \cite{G07}.

By the way, there exists another contribution by N.~N.~Nekhoroshev to Hamiltonian dynamics that does admit a reversible counterpart. It is his theorem on smooth $k$-parameter families of isotropic invariant $k$-tori carrying conditionally periodic motions in Hamiltonian systems with $n\geq k$ degrees of freedom and $k$ independent integrals in involution \cite{N94} (see also \cite{BG02,G02} for detailed proofs). D.~Bambusi \cite{B15} has carried over Nekhoroshev's result of \cite{N94} to reversible systems.

One may suggest that the Hamiltonian analogue of the distinction between the reversible contexts~1 and~2 is the ``splitting'' of the Hamiltonian KAM theory into the context of isotropic invariant tori and several contexts of non-isotropic invariant tori, see \cite{BS10,BHS96Gro,BHS96LNM,S08,KP01,Y92} and references therein. However, such an analogy seems very artificial and far-fetched. For instance, the persistence of invariant tori in the reversible context~2 requires the presence of external parameters. This is, generally speaking, not the case for coisotropic or atropic (neither isotropic nor coisotropic) invariant tori in Hamiltonian systems.

It is hardly feasible to unify the Hamiltonian KAM theory and the reversible one as a whole (cf.\ Problem~10 in \cite{S08}) although such a unification is probably achievable for the Hamiltonian isotropic context and the reversible context~1.

\section*{Acknowledgments}

It is the author's pleasure to thank H.~W.~Broer and his former students H.~Han{\ss}mann and G.~B.~Huitema for a long-term friendship and many fruitful discussions on various aspects of KAM theory. The study was partially supported by a grant of the President of the Russian Federation, project No.\ NSh-9789.2016.1.

\end{document}